\documentclass[12pt]{amsart}
\usepackage{latexsym,amssymb,amsrefs,graphicx}
\usepackage[all]{xy}

\newtheorem{lemma}{Lemma}[section]

\newtheorem{proposition}[lemma]{Proposition}
\newtheorem{theorem}[lemma]{Theorem}
\newtheorem{corollary}[lemma]{Corollary}
\newtheorem{definition}[lemma]{Definition}

\theoremstyle{definition}

\newcommand{\textmc}[1]{\textsf{#1}}

\begin{document}

\title[Three approaches to Morse-Bott homology]{Three approaches to Morse-Bott homology}

\author{David E. Hurtubise}
\address{Department of Mathematics and Statistics\\
         Penn State Altoona\\
         Altoona, PA 16601-3760}
\email{Hurtubise@psu.edu}

\dedicatory{Dedicated to Professor Augustin Banyaga\\
            on the occasion of his $65$th birthday}

\subjclass[2010]{Primary: 57R70 Secondary: 37D15 58E05}

\begin{abstract}
In this paper we survey three approaches to computing the homology of a 
finite dimensional compact smooth closed manifold using a Morse-Bott
function and discuss relationships among the three approaches.
The first approach is to perturb the function to a Morse function,
the second approach is to use moduli spaces of cascades, and the
third approach is to use the Morse-Bott multicomplex. 

With respect to an explicit perturbation (which can be used to derived
the Morse-Bott inequalities), the first two approaches yield the same
chain complex up to sign. The third approach is fundamentally different.
It combines singular cubical chains and Morse chains in the same multicomplex,
which provides a way of interpolating between the singular cubical chain
complex and the Morse-Smale-Witten chain complex.
\end{abstract}

\maketitle


\section{Introduction}

Functions with non-isolated critical points often arise in situations 
where there is some sort of symmetry or a group action. In such
situations the concept of a Morse-Bott function and the homology groups
associated to a Morse-Bott function developed by Raoul Bott in the 1950's
have proved to be extremely useful \cite{BotApp} \cite{BotNon} \cite{BotAna}.

\smallskip
A Morse-Bott function $f:M \rightarrow \mathbb{R}$ on a finite dimensional
compact smooth closed manifold $M$ is a function that is nondegenerate in the
directions normal to its critical submanifolds (Definition \ref{MorseBott}).
Bott found several applications of relationships he discovered between the
Betti numbers of the critical submanifolds of the Morse-Bott function $f$ and
the Betti numbers of the underlying manifold $M$. These relationships are
encoded in the Morse-Bott inequalities (Theorem \ref{MorseBottineq}),
which generalize the Morse inequalities (cf. Section 3.4 of \cite{BanLec}).
A Morse function is a Morse-Bott function with isolated critical points, and
the Morse inequalities give relationships between the critical points
of the Morse function and the Betti numbers of the underlying manifold.

The theory developed by Marston Morse in the 1920's originally gave information about the
Betti numbers of manifolds \cite{MorRel}, but more modern versions of Morse theory
determine CW-structures and chain complexes. In particular, the Morse Homology Theorem, which
was proved several decades after the Morse inequalities, says that the critical points and
gradient flow lines of a Morse function determine a chain complex that computes the homology
of the underlying manifold. The Morse inequalities follow as a direct consequence of the Morse
Homology Theorem (cf. Section 3 of \cite{BanDyn}), and thus one would expect that the
chain complex determined by a Morse function should contain more information
and have more applications than the Morse inequalities.  This is indeed the case.
Numerous applications of the Morse chain complex and its infinite dimensional analogues,
pioneered by Floer, have been found in many different contexts.

Correspondingly, one would expect that a chain complex determined by a Morse-Bott
function that computes the homology of the underlying manifold should contain
more information and have more applications than the Morse-Bott inequalities.
However, to find these applications one must first understand how a Morse-Bott
function determines a chain complex. This understanding is complicated by the fact
that there are several different chain complexes associated to a Morse-Bott function,
and there are multiple ways that a Morse-Bott function can be used to compute the
homology of the underlying manifold.

\smallskip
In this paper we will survey three of the many approaches that have been discovered
for using a Morse-Bott function to construct a chain complex that computes the homology
of the underlying manifold. The following three approaches to Morse-Bott homology will
be discussed, together with relationships among the three approaches.
\begin{enumerate}
\item Perturb the Morse-Bott function to a Morse function and use the
      Morse-Smale-Witten chain complex of the perturbed function.
\item Pick auxiliary Morse functions on the critical submanifolds and use
      the cascade chain complex defined using the auxiliary functions.
\item Use the Morse-Bott multicomplex, which provides a common framework for both
      Morse chains and singular cubical chains.
\end{enumerate}

In Section \ref{perturbsection} we will explain how the first approach can be used
to prove the Morse-Bott inequalities following \cite{BanDyn}. In Section
\ref{cascadesection} we will describe the construction of the cascade chain complex
and explain how the cascade chain complex can be viewed as the Morse-Smale-Witten
chain complex of a specific perturbation of the Morse-Bott function following
\cite{BanCas}. Section \ref{multicomplexsection} will discuss the Morse-Bott
multicomplex developed by Banyaga and Hurtubise using singular cubical chains
and fibered product constructions \cite{BanMor}. The Morse-Bott multicomplex
is fundamentally different from other approaches to Morse-Bott
homology. It provides a common framework for singular cubical chains
and Morse chains, and thus it provides a way of interpolating between 
the singular cubical chain complex and the Morse-Smale-Witten chain complex.

\smallskip
Other approaches to Morse-Bott homology that will not be discussed in detail in this
paper include: the spectral sequence associated with the filtration of the manifold
determined by a Morse-Bott function, the de Rham cochain complex of Austin and Braam 
\cite{AusMor}, and the Morse-Bott chain complex defined using currents due to Latschev 
\cite{LatGra}. \nocite{ChoMor} One common distinguishing feature of both Morse and 
Morse-Bott homology is that the boundary operators that define the homology are expressed
in terms of moduli spaces of gradient flow lines. We will not discuss the spectral sequence
associated with the filtration of the manifold because currently there is no known way to
relate the differentials in that spectral sequence to moduli spaces of gradient flow lines.

The paper by Austin and Braam \cite{AusMor} was a source of inspiration for the results
discussed in Sections \ref{perturbsection} and \ref{multicomplexsection}, and connections
with their work are discussed in those sections. (A paper by Fukaya \cite{FukFlo} also 
served as a source of inspiration for the results discussed in Section
\ref{multicomplexsection}.) However, the Austin-Braam approach uses differential forms
to construct a comulticomplex that computes the de Rham cohomology of the manifold 
with real coefficients, whereas the three approaches discussed in detail in this 
paper all concern homology with integer coefficients. Finally, the Morse-Bott chain complex
defined by Latschev \cite{LatGra} is part of the program of Harvey and Lawson \cite{HarFin}
to approach Morse theory using the de Rham-Federer theory of currents. Although their
approach to Morse theory is very interesting, it is also quite different from the other
approaches discussed in this paper and will not be further reviewed.


\section{The Morse-Smale-Witten chain complex}\label{MorseChain}

In this section we briefly recall the construction of the Morse-Smale-Witten
chain complex and the Morse Homology Theorem.  For more details see \cite{BanLec}.

\smallskip
Let $Cr(f) = \{p \in M |\, df_p = 0 \}$ denote the set of critical points of a 
smooth function $f:M \rightarrow \mathbb{R}$ on a finite dimensional smooth 
manifold $M$. A critical point $p \in Cr(f)$ is said to be \textbf{nondegenerate} if
and only if the Hessian $H_p(f)$ is nondegenerate. The \textbf{index} $\lambda_p$
of a nondegenerate critical point $p$ is defined to be the index of the symmetric
bilinear form $H_p(f)$, i.e. the dimension of the subspace of $T_pM$
where $H_p(f)$ is negative definite. If all the critical points of $f$ are
non-degenerate, then $f$ is called a \textbf{Morse function}.

If $f:M \rightarrow \mathbb{R}$ is a Morse function on an $m$-dimensional 
compact smooth Riemannian manifold $(M,g)$, then the \textbf{stable manifold} $W^s_f(p)$
and the \textbf{unstable manifold} $W^u_f(p)$ of a critical point $p \in Cr(f)$ are
defined to be
\begin{eqnarray*}
W^s_f(p) & = & \{ x\in M | \lim_{t \rightarrow \infty} \varphi_t(x) = p \}\\
W^u_f(p) & = & \{ x\in M | \lim_{t \rightarrow -\infty} \varphi_t(x) = p \}
\end{eqnarray*}
where $\varphi_t$ is the 1-parameter group of diffeomorphisms generated by
minus the gradient vector field, i.e. $-\nabla f$. The Stable/Unstable Manifold Theorem
for a Morse Function says that the tangent space at $p$ splits as 
$$
T_pM = T^s_pM \oplus T_p^uM
$$
where the Hessian is positive definite on $T_p^sM \stackrel{def}{=} T_p W^s_f(p)$ and 
negative definite on $T_p^uM \stackrel{def}{=} T_p W^u_f(p)$.  Moreover, the
stable and unstable manifolds of $p$ are surjective images of smooth 
embeddings
\begin{eqnarray*}
E^s: T_p^sM & \rightarrow & W^s_f(p) \subseteq M\\
E^u: T_p^uM & \rightarrow & W^u_f(p) \subseteq M. 
\end{eqnarray*}
Hence, $W^s_f(p)$ is a smoothly embedded open disk of dimension $m - \lambda_p$, 
and $W^u_f(p)$ is a smoothly embedded open disk of dimension $\lambda_p$.

If the stable and unstable manifolds of a Morse function $f:M \rightarrow
\mathbb{R}$ all intersect transversally, then the function $f$ is called
\textbf{Morse-Smale}. For any metric $g$ on $M$ the set of smooth Morse-Smale
functions is dense in the space of all smooth functions on $M$ by the Kupka-Smale 
Theorem (cf. Theorem 6.6 and Remark 6.7 of \cite{BanLec}), and for a given Morse function 
$f:M \rightarrow \mathbb{R}$ one can always find a Riemannian metric on $M$ so that $f$
is Morse-Smale with respect to the chosen metric (cf. Theorem 2.20 of \cite{AbbLec}).
Moreover, if $f$ is Morse-Smale and $p,q\in Cr(f)$ then $W_f(q,p) = W^u_f(q) \cap W^s_f(p)$
is an embedded submanifold of $M$ of dimension $\lambda_q  - \lambda_p$, and when
$\lambda_q - \lambda_p = 1$ the number of gradient flow lines from $q$ to $p$ is finite
(cf. Corollary 6.29 of \cite{BanLec}).

If we choose an orientation for each of the unstable manifolds of $f$, then there is an
induced orientation on the normal bundles of the stable manifolds. Thus, we can define 
an integer associated to any two critical points $p$ and $q$ of relative index one by 
counting the number of gradient flow lines from $q$ to $p$ with signs determined by the
orientations. This integer is denoted by $n_f(q,p) = \# \mathcal{M}_f(q,p)$, where 
$\mathcal{M}_f(q,p) = W_f(q,p)/\mathbb{R}$ is the moduli space of gradient flow lines of
$f$ from $q$ to $p$. The \textbf{Morse-Smale-Witten chain complex} is defined to be the
chain complex $(C_\ast(f),\partial_\ast)$ where $C_k(f)$ is the free abelian group
generated by the critical points $q$ of index $k$ and the boundary operator 
$\partial_k:C_k(f) \rightarrow C_{k-1}(f)$ is given by
$$
\partial_k(q)\ \ = \sum_{p \in Cr_{k-1}(f)} n_f(q,p)p
$$
where $Cr_{k-1}(f)$ denotes the set of critical points with index $k-1$.

\begin{theorem}[Morse Homology Theorem]\label{Morsehomology}
The pair $(C_\ast(f),\partial_\ast)$ is a chain complex, and the homology
of $(C_\ast(f),\partial_\ast)$ is isomorphic to the singular homology 
$H_\ast(M;\mathbb{Z})$.
\end{theorem}

\noindent
Note that the Morse Homology Theorem implies that the homology of
$(C_\ast(f),\partial_\ast)$ is independent of the Morse-Smale
function $f:M \rightarrow \mathbb{R}$, the Riemannian metric,
and the chosen orientations.

There are many different ways to prove the Morse Homology Theorem. The approach in 
\cite{BanLec} uses the Conley index and Conley's connection matrix to give an 
explicit isomorphism between the Morse homology groups and the singular homology groups.
Another approach is to first show that the unstable manifolds of a Morse-Smale
function $f:M \rightarrow \mathbb{R}$ determine a CW-structure $X$ on $M$, and then
show that the chain complex $(\underline{C}_\ast(X),\partial^{CW}_\ast)$ determined by the
CW-structure is the same as the Morse-Smale-Witten chain complex $(C_\ast(f),\partial_\ast)$.
Both of these steps are nontrivial. For instance, $\partial^{CW}_\ast$ is induced from a
connecting homomorphism in the homology exact sequence of a triple, whereas $\partial_\ast$
is defined by counting gradient flow lines. (See Section 7.1 of \cite{BanLec} for a
more complete discussion of the technical details encountered when using this approach.)

However, in spite of the difficulties, it is possible to prove that the unstable
manifolds of a Morse-Smale function $f:M \rightarrow \mathbb{R}$ determine a CW-structure
$X$ on $M$ and the following diagram commutes
$$
\xymatrix{
C_k(f) \ar@{<->}[d]^{\approx} \ar[r]^{\partial_k} & C_{k-1}(f) \ar@{<->}[d]^{\approx}\\
\underline{C}_k(X) \ar[r]^{\partial^{CW}_k}       & \underline{C}_{k-1}(X)
}
$$
for every $k \in \mathbb{Z}_+$, where the vertical maps are induced by identifying
critical points with their unstable manifolds. For more details concerning this approach
see the recent papers by Qin \cite{QinOnm} \cite{QinAnA} and the references therein.
For a list of other approaches to proving the Morse Homology Theorem see the introduction
to \cite{BanLec}.


\section{Perturbing a Morse-Bott function to a Morse function}\label{perturbsection}

The chain groups in the Morse-Smale-Witten chain complex are finitely
generated because a Morse function $f:M \rightarrow \mathbb{R}$ on a finite
dimensional compact manifold has a finite number of isolated critical points.
If the critical points of $f$ are not isolated, then $f$ has an infinite
number of critical points and the Morse Homology Theorem does not apply.  
In this case, some additional assumptions and/or auxiliary data are required to
construct a chain complex, a multicomplex, or a spectral sequence.

\medskip
Let $f:M \rightarrow \mathbb{R}$ be a smooth function whose critical set
$\text{Cr}(f)$ contains a submanifold $C$ of positive dimension. 
Pick a Riemannian metric on $M$ and use it to split $T_\ast M|_C$
as
$$
T_\ast M|_C = T_\ast C \oplus \nu_\ast C
$$
where $T_\ast C$ is the tangent space of $C$ and $\nu_\ast C$ is the normal bundle of $C$.
Let $p \in C$, $V \in T_p C$, $W \in T_pM$, and let $H_p(f)$ be the Hessian of $f$ at $p$.
We have
$$
H_p(f)(V,W) = V_p \cdot (\tilde{W} \cdot f) = 0
$$
since $V_p \in T_pC$ and any extension of $W$ to a vector field $\tilde{W}$ satisfies 
$df(\tilde{W})|_C$ $= 0$.  Therefore, the Hessian $H_p(f)$ induces a
symmetric bilinear form $H_p^\nu(f)$ on $\nu_p C$.

\begin{definition}\label{MorseBott}
A smooth function $f:M \rightarrow \mathbb{R}$ on a smooth manifold $M$ is 
called a \textbf{Morse-Bott function} if and only if the set of critical points
$\text{\rm Cr}(f)$ is a disjoint union of connected submanifolds and for each connected
submanifold $C \subseteq \text{\rm Cr}(f)$ the bilinear form $H_p^\nu(f)$ is non-degenerate 
for all $p \in C$. The Morse-Bott \textbf{index} of a critical submanifold 
$C \subseteq \text{\rm Cr}(f)$ is defined to be the index of $H_p^\nu(f)$ for any $p \in C$.
\end{definition}

\noindent
Note:  The Morse-Bott index is well defined by the Morse-Bott Lemma (cf. Section 3.5 of
\cite{BanLec}).

\smallskip
A well known theorem says that on a compact closed smooth manifold the space of smooth
Morse functions is open and dense in the space of all smooth functions (cf. Theorem 5.31
of \cite{BanLec}). So, one approach to computing homology from a Morse-Bott function would
be to perturb the Morse-Bott function to a Morse function and apply Theorem
\ref{Morsehomology} (Morse Homology Theorem) using a metric on $M$ such that the perturbed 
function is Morse-Smale with respect to the chosen metric.

\smallskip
Perturbing a Morse-Bott function to a Morse function using abstract perturbations and
defining a Morse-like chain complex associated to the perturbed function is a standard
technique used in gauge theory with respect to the Chern-Simons functional and in Floer
theory with respect to the symplectic action functional. In the setting of a Morse-Bott
function $f:M \rightarrow \mathbb{R}$ on a finite dimensional orientable compact smooth
manifold $M$, Austin and Braam defined a more explicit perturbation of 
$f:M \rightarrow \mathbb{R}$ to a Morse function $h:M \rightarrow \mathbb{R}$ by 
introducing Morse functions on the critical submanifolds \cite{AusMor}.

Austin and Braam used their perturbation technique to compare the homology of a filtered
cochain complex $(C^\ast,\partial^\ast)$ they defined using differential forms on the critical
submanifolds of a Morse-Bott-Smale function $f:M \rightarrow \mathbb{R}$ (Definition
\ref{MBStransversality}) with the Morse-Smale-Witten cochain complex
$(C^\ast(h) \otimes \mathbb{R}, \partial^\ast)$ of the perturbed function 
$h:M \rightarrow \mathbb{R}$ (Proposition 3.10 of \cite{AusMor}). Austin and Braam's 
cochain complex $(C^\ast,\partial^\ast)$ has the structure of a comulticomplex 
(Definition \ref{multicomplex}), which determines a spectral sequence coming from the filtration 
\cite{HurMul}. By exhibiting a chain morphism between the filtered cochain complexes
$(C^\ast,\partial^\ast)$ and $(C^\ast(h) \otimes \mathbb{R}, \partial^\ast)$ that
induces an isomorphism of the $E^1$ terms of the spectral sequences determined by the
filtrations, Austin and Braam prove that there is also an isomorphism on the $E^\infty$
terms, and hence an induced isomorphism on the homology of the filtered cochain complexes.
This proves that both cochain complexes compute the de Rham cohomology of $M$.

Corollary 3.9 of \cite{AusMor} states without proof that the polynomial Morse-Bott 
inequalities follow from the fact that the comulticomplex $(C^\ast,\partial^\ast)$
computes the de Rham cohomology of $M$. A proof of the Morse-Bott inequalities along those
lines would most likely involve an analysis of the spectral sequence determined by the
comulticomplex, and hence would not be as immediate as proving that the polynomial
Morse inequalities follow from the existence of a CW-complex determined by the Morse function
or the Morse Homology Theorem (cf. Section 3.4 of \cite{BanLec} and Section 3 of
\cite{BanDyn}). However, Banyaga and Hurtubise showed in \cite{BanDyn} that it is possible
to apply the perturbation technique used by Austin and Braam together with the polynomial
Morse inequalities to prove the polynomial Morse-Bott inequalities without appealing to the
full Morse-Bott multicomplex.


\subsection*{The polynomial Morse-Bott inequalities}
Let $M$ be a compact smooth manifold of dimension $m$, and define the $k^\text{th}$
\textbf{Betti number} of $M$, denoted $b_k$, to be the rank of the $k^\text{th}$ homology
group $H_k(M;\mathbb{Z})$ modulo its torsion subgroup. Let $f:M \rightarrow \mathbb{R}$ 
be a smooth Morse function on $M$, and let $\nu_k$ denote the number of critical points of
$f$ of index $k$ for all $k=0,\ldots , m$.

\begin{definition}
The \textbf{Poincar\'e polynomial} of $M$ is defined to be
$$
P_t(M) = \sum_{k=0}^m b_k t^k,
$$
and the \textbf{Morse polynomial} of $f$ is defined to be
$$
M_t(f) = \sum_{k=0}^m \nu_k t^k.
$$
\end{definition}

For a proof of the following theorem using the fact that $f:M \rightarrow \mathbb{R}$
determines a CW-complex homotopic to $M$ with $\nu_k$ cells of dimension
$k$ for all $k=0, \ldots , m$ see Section 3.4 of \cite{BanLec}.  For a similar proof
that uses the Morse Homology Theorem instead of the CW-complex see Section 3 of
\cite{BanMor}.

\begin{theorem}[Polynomial Morse Inequalities]\label{Morsepolynomial}
For any Morse function $f:M \rightarrow \mathbb{R}$ on a compact smooth manifold
of dimension $m$ we have
$$
M_t(f) = P_t(M) + (1+t) R(t)
$$
where $R(t)$ is a polynomial with non-negative integer coefficients.
That is, $R(t) = \sum_{k=0}^{m-1} r_k t^k$ where $r_k \in \mathbb{Z}$ 
satisfies $r_k \geq 0$ for all $k=0,\ldots ,m-1$.
\end{theorem}

\smallskip
Now let $f:M \rightarrow \mathbb{R}$ be a Morse-Bott function, and assume that
$$
\text{Cr}(f) = \coprod_{j=1}^l C_j,
$$
where $C_1,\ldots ,C_l$ are disjoint connected critical submanifolds.  

\begin{definition}
The Morse-Bott polynomial of $f$ is defined to be
$$
MB_t(f) = \sum_{j=1}^l P_t(C_j) t^{\lambda_j}
$$
where $\lambda_j$ is the Morse-Bott index of the critical submanifold $C_j$
and $P_t(C_j)$ is the Poincar\'e polynomial of $C_j$.
\end{definition}

\smallskip
Bott proved a version of the following result stated in terms of Betti numbers of
homology with local coefficients in an orientation bundle in place of any
orientation assumptions \cite{BotNon} \cite{BotLec}.

\begin{theorem}[Morse-Bott Inequalities]\label{MorseBottineq}
Let $f:M \rightarrow \mathbb{R}$ be a Morse-Bott function on a finite
dimensional orientable compact smooth manifold, and assume that all the critical 
submanifolds of $f$ are orientable. Then there exists a polynomial 
$R(t)$ with non-negative integer coefficients such that
$$
MB_t(f) = P_t(M) + (1+t)R(t).
$$
\end{theorem}

\smallskip
Bott proved his version of this theorem by studying how the homotopy type of the
``half-spaces'' $M^y = f^{-1}(-\infty ,y]$ change as $y$ crosses critical values. 
In \cite{BanDyn}, Banyaga and Hurtubise gave a proof of Theorem \ref{MorseBottineq}
via a dynamical systems approach by expanding on the perturbation technique used by 
Austin and Braam in \cite{AusMor}.


\subsection*{Outline of the Banyaga-Hurtubise proof}
Chose a small tubular neighborhood $T_j$ around each connected component
$C_j\subseteq \text{Cr}(f)$ for all $j=1,\ldots ,l$ with local coordinates $(u,v,w)$ 
consistent with those from the Morse-Bott Lemma (cf. Section 3.5 of \cite{BanLec} 
or \cite{BanApr}). Pick a Riemannian metric on $M$ such that the charts from the
Morse-Bott Lemma are isometries with respect to the standard Euclidean metric on
$\mathbb{R}^m$, and then pick positive Morse functions $f_j:C_j \rightarrow \mathbb{R}$
that are Morse-Smale with respect to the restriction of the Riemannian metric to $C_j$
for all $j=1,\ldots ,l$. The Morse-Smale functions $f_j:C_j \rightarrow \mathbb{R}$
exist by the Kupka-Smale Theorem.

For every $j=1,\ldots ,l$ extend $f_j$ to a function on $T_j$ by making $f_j$ constant
in the direction normal to $C_j$. Let $\tilde{T}_j \subset T_j$ be a smaller tubular
neighborhood of $C_j$ with the same coordinates as $T_j$, and let $\rho_j$ 
be a smooth nonincreasing bump function that is constant in the direction parallel
to $C_j$, equal to $1$ on $\tilde{T}_j$, and equal to $0$ outside of $T_j$.
For a small $\varepsilon > 0$ the function $h:M \rightarrow \mathbb{R}$ given by 
$$
h = f + \varepsilon \left( \sum_{j=1}^l \rho_j f_j \right)
$$
is a Morse function close to $f$, and the critical points of $h$ are exactly the critical
points of the $f_j$ for $j=1,\ldots ,l$. Moreover, if $p \in C_j$ is a critical point
of $f_j:C_j \rightarrow \mathbb{R}$ of index $\lambda_p^j$, then $p$ is a 
critical point of $h$ of index $\lambda_p^h = \lambda_j + \lambda_p^j$.
A well-known folk theorem (cf. Section 2.12 of \cite{AbbLec}) says that it is possible
to perturb the Riemannian metric on $M$ outside of the union of the tubular neighborhoods 
$T_j$ for $j=1,\ldots ,l$ so that $h$ satisfies the Morse-Smale transversality condition 
with respect to the perturbed metric. 

This explicit perturbation and choice of metric makes it possible to compare the
Morse-Smale-Witten chain complex of $h$ with those of $f_j$ for $j=1,\ldots, l$. In 
particular, for every $n=0,\ldots ,m$ we have the following description of the $n^{\text{th}}$
Morse-Smale-Witten chain group of $h$ in terms of the Morse-Smale-Witten chain groups of
the $f_j$ for $j=1,\ldots ,l$.
$$
C_n(h) = \bigoplus_{\lambda_j + k = n} C_k(f_j)
$$
Now let $M_t(f_j)$ denote the Morse polynomial of $f_j:C_j\rightarrow \mathbb{R}$, and note
that the relation $\lambda_p^h = \lambda_j+ \lambda_p^j$ implies that
$$
M_t(h) = \sum_{j=1}^l M_t(f_j) t^{\lambda_j}.
$$
The polynomial Morse inequalities (Theorem \ref{Morsepolynomial}) say that
$$
M_t(h) = P_t(M) + (1+t)R_h(t)
$$
and
$$
M_t(f_j) = P_t(C_j) + (1+t)R_j(t)
$$
where $R_h(t)$ and $R_j(t)$ are polynomials with non-negative integer
coefficients for all $j=1,\ldots ,l$.  This leads to the following
straightforward computation.
\begin{eqnarray*}
MB_t(f) & =  & \sum_{j=1}^l P_t(C_j) t^{\lambda_j}\\
& = & \left.\left.\sum_{j=1}^l \right( M_t(f_j) - (1+t) R_j(t) \right)t^{\lambda_j}\\
& = & \sum_{j=1}^l M_t(f_j)t^{\lambda_j} - (1+t) \sum_{j=1}^l R_j(t)t^{\lambda_j}\\
& = & M_t(h) - (1+t) \sum_{j=1}^l R_j(t)t^{\lambda_j}\\
& = & P_t(M) + (1+t) R_h(t) - (1+t) \sum_{j=1}^l R_j(t)t^{\lambda_j}\\
& = & P_t(M) + (1+t) \left( R_h(t) - \sum_{j=1}^l R_j(t)t^{\lambda_j} \right)
\end{eqnarray*}

It remains to show that the polynomial multiplying $(1+t)$ in the last line has 
non-negative integer coefficients. This is accomplished by first noting that the proof
of the polynomial Morse inequalities from the Morse Homology Theorem shows that the 
polynomial $R_j(t)$ is given by 
$$
R_j(t) = \sum_{k=1}^{c_j} (\nu_k^j - z_k^j)t^{k-1}
$$
where $c_j = \text{dim } C_j$, $\nu_k^j = \text{rank }C_k(f_j)$, and $z_k^j$ is the rank of
the kernel of the boundary operator $\partial_k^{f_j}:C_k(f_j) \rightarrow C_{k-1}(f_j)$
in the Morse-Smale-Witten chain complex of $f_j:C_j\rightarrow \mathbb{R}$.  The proof 
of Theorem \ref{MorseBottineq} is then completed by using the close relationship between
the dynamics of the gradient flow lines of $f$, $f_j$, and $h$ to show that 
$\sum_{\lambda_j + k = n} z_k^j \geq z_n^h$ for all $n=1,\ldots ,m$, where $z_n^h$
is the rank of the kernel of the boundary operator $\partial_n^h:C_k(h) \rightarrow
C_{k-1}(h)$ in the Morse-Smale-Witten chain complex of $h:M \rightarrow \mathbb{R}$.


\subsection*{Comparison with Bott's proof}
Bott's proof of the polynomial Morse-Bott inequalities was based on studying what he
called ``half-spaces'' $M^y = f^{-1}(-\infty, y]$ \cite{BotMor}. The homotopy type of 
$M^y$ is the same on any interval $a < y < b$ that doesn't contain a critical value, 
and when $y$ crosses a critical value $c$ the homotopy type changes by the attachment of
disk bundles whose dimensions are given by the Morse-Bott indexes of the critical
submanifolds in the level set $f^{-1}(c)$ (cf. Appendix B of \cite{FarTop}).

Bott's original version of the polynomial Morse-Bott inequalities avoided any orientation
assumptions by using Betti numbers with local coefficients in an orientation bundle in place
of the standard Betti numbers. The orientation bundles of the critical submanifolds 
he considered are determined by disk bundles given by the unstable part of the gradient
flow near the critical submanifolds. When these disk bundles are orientable the Betti
numbers with local coefficients in the orientation bundles reduce to the standard Betti
numbers of the critical submanifolds. For more details see the original papers by Bott or 
Appendix C of \cite{FarTop}.

It is interesting to note that Bott's version of the Morse-Bott inequalities reduces to
the conclusion of Theorem \ref{MorseBottineq} when the disk bundles given by the
negative part of the gradient flow near the submanifolds are orientable, whereas the
proof given by Banyaga and Hurtubise assumes that $M$ and the critical submanifolds are
orientable. That is, the tangent space of $M$ along a critical submanifold $C$
has a decomposition
$$
T_\ast M = T_\ast C \oplus \nu_\ast^-C \oplus \nu_\ast^+C,
$$
and Banyaga and Hurtubise assumed that $T_\ast M$ and $T_\ast C$ are orientable in
order to prove Theorem \ref{MorseBottineq}.  On the other hand, it is the assumption
that the bundle $\nu_\ast^-C$ is orientable that allows one to conclude that
the Betti numbers with local coefficients in the orientation bundle used by Bott reduce
to the Betti numbers considered in Theorem \ref{MorseBottineq}. These two conditions
are distinct when $\nu_\ast^+C$ is not orientable.


\section{Cascades}\label{cascadesection}

A second approach to computing homology using a Morse-Bott function involves
introducing Morse functions on the critical submanifolds and defining chain
groups generated by the critical points of the Morse functions that agree with
those defined in the previous section. However, the boundary operator is defined
by counting the number of ``cascades'' between two critical points of relative
index one, which are defined without reference to the perturbed function. Roughly
speaking, a cascade between two critical points is a concatenation of some gradient
flow lines of the Morse-Bott function and pieces of the gradient flow lines of the Morse 
functions on the critical submanifolds. Cascades were introduced by Frauenfelder in
\cite{FraTheA} \cite{FraThesis}, and cascade-like objects were introduced independently
in the context of holomorphic curves by Bourgeois in \cite{BouThesis} \cite{BouAMo}.
Cascades have since been used by several authors studying symplectic and contact homology
\cite{BouAne} \cite{BouSym} \cite{BouAMo} \cite{CieAFl}. We begin our discussion of cascades
with the following definitions from \cite{BanCas}.


\subsection*{Cascades and Morse-Bott-Smale transversality}
Let $f:M \rightarrow \mathbb{R}$ be a Morse-Bott function on a finite dimensional
compact smooth manifold, and let
$$
\text{Cr}(f) = \coprod_{j=1}^l C_j,
$$
where $C_1,\ldots ,C_l$ are disjoint connected critical submanifolds of Morse-Bott index
$\lambda_1,\ldots ,\lambda_l$ respectively. Let $f_j:C_j \rightarrow \mathbb{R}$ be 
a Morse function on the critical submanifold $C_j$ for all $j=1,\ldots ,l$. For 
a critical point $q \in C_j$ of $f_j:C_j\rightarrow \mathbb{R}$ denote the Morse index of
$q$ relative to $f_j$ by $\lambda_q^j$, the stable manifold of $q$ relative to $f_j$
by $W^s_{f_j}(q) \subseteq C_j$, and the unstable manifold of $q$ relative to $f_j$ by
$W^u_{f_j}(q) \subseteq C_j$.

\begin{definition}\label{index}
If $q\in C_j$ is a critical point of the Morse function $f_j:C_j \rightarrow \mathbb{R}$
for some $j=1,\ldots ,l$, then the \textbf{total index} of $q$, denoted $\lambda_q$, is
defined to be the sum of the Morse-Bott index of $C_j$ and the Morse index of $q$ relative
to $f_j$, i.e.
$$
\lambda_q = \lambda_j + \lambda_q^j.
$$
\end{definition}

\begin{definition}\label{flowlinecascade}
For $q\in \text{Cr}(f_j)$, $p \in \text{Cr}(f_i)$, and $n \in \mathbb{N}$, a 
\textbf{flow line with $n$ cascades from $q$ to $p$} is a $2n-1$-tuple:
$$
\left((x_k)_{1\leq k\leq n},(t_k)_{1 \leq k \leq n-1}  \right)
$$
where $x_k \in C^\infty(\mathbb{R},M)$ and $t_k \in \mathbb{R}_+ =
\{t \in \mathbb{R} |\ t \geq 0\}$ satisfy the following for all $k$.
\begin{enumerate}
\item Each $x_k$ is a non-constant gradient flow line of $f$, i.e.
$$
\frac{d}{dt}x_k(t) = - (\nabla f)(x_k(t)).
$$
\item For the first cascade $x_1(t)$ we have 
$$
\lim_{t \rightarrow -\infty} x_1(t) \in W^u_{f_j} (q) \subseteq C_j,
$$
and for the last cascade $x_n(t)$ we have 
$$
\lim_{t \rightarrow \infty} x_n(t) \in W^s_{f_i}(p) \subseteq C_i.
$$
\item For $1 \leq k \leq n-1$ there are critical submanifolds $C_{j_k}$ and
gradient flow lines $y_k \in C^\infty(\mathbb{R},C_{j_k})$ of $f_{j_k}$, i.e. 
$$
\frac{d}{dt}y_k(t) = - (\nabla f_{j_k})(y_k(t)),
$$
such that $\lim_{t \rightarrow \infty} x_k(t) = y_k(0)$ and 
$\lim_{t \rightarrow -\infty} x_{k+1}(t) = y_{k}(t_k)$.
\end{enumerate}
When $j = i$ a \textbf{flow line with zero cascades from $q$ to $p$} is a gradient flow
line of $f_j$ from $q$ to $p$.
\end{definition}

\begin{figure}[h]
\includegraphics{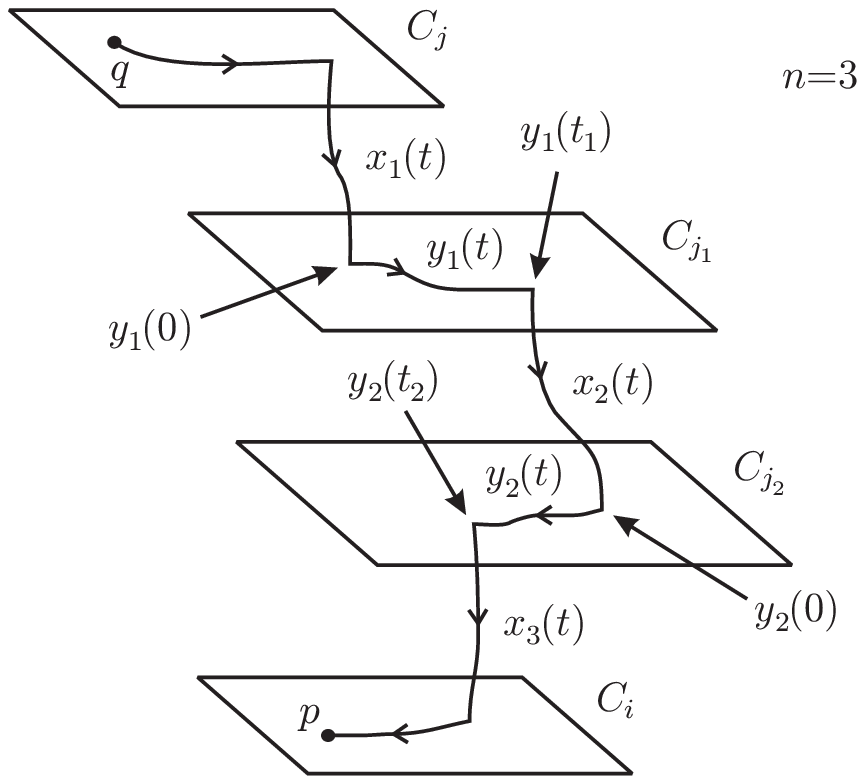}
\end{figure}

\smallskip\noindent
Note: In the preceding definition the parameterizations of the gradient
flow lines $y_k(t)$ of the Morse functions $f_{j_k}:C_{j_k}\rightarrow \mathbb{R}$ are
fixed in $(3)$ by $\lim_{t \rightarrow \infty} x_k(t) = y_k(0)$,
and the entry $t_k$ records the time spent flowing along the critical submanifold
$C_{j_k}$ (or resting at a critical point). However, the parameterizations of the 
cascades $x_1(t),\ldots ,x_n(t)$ are not fixed. Hence, there is an action of
$\mathbb{R}^n$ on a flow line with $n$ cascades given by
$$
\left((x_k(t))_{1\leq k\leq n},(t_k)_{1 \leq k \leq n-1}  \right) \mapsto
\left((x_k(t+s_k))_{1\leq k\leq n},(t_k)_{1 \leq k \leq n-1}  \right)
$$
for $(s_1,\ldots ,s_n) \in \mathbb{R}^n$.

\begin{definition}
For $q\in \text{Cr}(f_j)$, $p \in \text{Cr}(f_i)$, and $n \in \mathbb{N}$ we denote
the space of flow lines from $q$ to $p$ with $n$ cascades by $W^c_n(q,p)$, and we 
denote the quotient of $W^c_n(q,p)$ by the action of $\mathbb{R}^n$ by
$$
\mathcal{M}^c_n(q,p) = W^c_n(q,p)/\mathbb{R}^n.
$$
The \textbf{set of unparameterized flow lines with cascades from $q$ to $p$} is defined
to be
$$
\mathcal{M}^c(q,p) = \bigcup_{n \in \mathbb{Z}_+} \mathcal{M}^c_n(q,p)
$$
where $\mathcal{M}^c_0(q,p) = W^c_0(q,p)/\mathbb{R}$. We will say 
that an element of $\mathcal{M}^c(q,p)$ \textbf{begins} at $q$ and \textbf{ends} at $p$.
\end{definition}

\smallskip
Somewhat surprisingly, under the right conditions moduli spaces of cascades 
have properties similar to moduli spaces of gradient flow lines of a Morse-Smale function.
In particular, under various assumptions it is possible to prove that $\mathcal{M}^c(q,p)$
is a smooth manifold of dimension $\lambda_q - \lambda_p - 1$, and every sequence in 
$\mathcal{M}^c(q,p)$ has a subsequence that converges to a broken flow line with cascades. 
These two fundamental properties imply that there are only a finite number of cascades
between any two critical points of relative index one, and hence it is possible to define
a boundary operator $\partial_\ast^c$ by counting cascades. (Proving directly that 
$\partial_\ast^c \circ \partial_\ast^c = 0$ requires a stronger result. Namely, that
$\mathcal{M}^c(q,p)$ has a compactification consisting of broken flow lines with cascades 
when $\lambda_q - \lambda_p = 2$.)

These fundamental properties were proved by Frauenfelder in \cite{FraTheA} under the
assumptions that the Riemannian metric on $M$ and the Riemannian metrics on the critical
submanifolds meet certain generic conditions that imply that a particular Fredholm operator
is surjective. They were also proved by Banyaga and Hurtubise in \cite{BanCas} under the
assumptions that the Morse-Bott function satisfies the Morse-Bott-Smale transversality
condition and the unstable and stable manifolds of the Morse functions on the
critical submanifolds are transverse to certain beginning and endpoint maps.

\begin{definition}[Morse-Bott-Smale Transversality]\label{MBStransversality}
A Morse-Bott function $f:M \rightarrow \mathbb{R}$ is said to satisfy the 
\textbf{Morse-Bott-Smale transversality} condition with respect to a given 
Riemannian metric $g$ on $M$ if and only if for any two connected critical submanifolds 
$C$ and $C'$, $W^u_f(q)$ intersects $W^s_f(C')$ transversely in $M$, i.e. $W^u_f(q)
\pitchfork W^s_f(C') \subseteq M$, for all $q\in C$.
\end{definition}

The Morse-Bott-Smale transversality condition depends on both the function and the
Riemannian metric, and it may not be possible to perturb the metric to make a 
given Morse-Bott function satisfy the Morse-Bott-Smale transversality condition.
For some interesting examples see Section 2 of \cite{LatGra}. This is quite different
from the situation encountered with the Morse-Smale transversality condition where 
it is always possible to perturb either the function or the metric to make the
condition hold. (Of course, one can always perturb the Morse-Bott function to
a Morse function to make the Morse-Bott-Smale transversality condition hold, but
that approach was already discussed in the previous section.)


\subsection*{Moduli spaces of gradient flow lines and cascades}
There are many technical consequences of the Morse-Bott-Smale transversality condition
that have implications for moduli spaces of gradient flow lines. For instance, the moduli
space of gradient flow lines between two critical submanifolds $C_k$ and $C_{k'}$ of a
Morse-Bott-Smale function $f:M \rightarrow \mathbb{R}$
$$
\mathcal{M}_f(C_k,C_{k'}) = \left(W^u_f(C_k) \cap W^s_f(C_{k'})\right)/\mathbb{R}
$$
is a manifold of dimension $\lambda_k - \lambda_{k'} + \text{dim }C_k - 1$
(cf. Lemma 3.5 of \cite{BanMor}), and the beginning point map
$$
\partial_-: \mathcal{M}_f(C_k,C_{k'}) \rightarrow C_k
$$
sending a gradient flow line to its starting point is a submersion (cf. Lemma 5.19 of
\cite{BanMor}). In fact, the moduli space of gradient flow lines $\mathcal{M}_f(C_k,C_{k'})$
has a compactification $\overline{\mathcal{M}}_f(C_k,C_{k'})$ consisting of broken
gradient flow lines, which is a smooth manifold with corners, and the beginning point map
$$
\partial_-: \overline{\mathcal{M}}_f(C_k,C_{k'}) \rightarrow C_k
$$
is both a submersion and a stratum submersion (cf. Corollary 5.20 of \cite{BanMor}).

These consequences of the Morse-Bott-Smale transversality condition were used by Banyaga
and Hurtubise to construct smooth manifolds with corners defined in terms of iterated
fibered products over the beginning and endpoint maps \cite{BanMor} and over the beginning
point map and the endpoint map composed with the gradient flow along the critical submanifolds
\cite{BanCas}. In the second case, the iterated fibered products can be viewed as spaces of
cascades from one critical submanifold to another. This leads to the condition that the
beginning and endpoint maps from the iterated fibered products are transverse to the unstable
and stable manifolds of the Morse functions on the critical submanifolds; a condition that is
always satisfied by an arbitrarily small perturbation of the Morse functions on the critical
submanifolds \cite{BanCas}.

\begin{theorem}\label{manifold}
Assume that $f$ satisfies the Morse-Bott-Smale transversality condition with respect to
the Riemannian metric $g$ on $M$, $f_k:C_k \rightarrow \mathbb{R}$ satisfies
the Morse-Smale transversality condition with respect to the restriction of $g$
to $C_k$ for all $k=1,\ldots ,l$, and the unstable and stable manifolds $W^u_{f_j}(q)$ and
$W^s_{f_{i}}(p)$ are transverse to the beginning and endpoint maps. 
\begin{enumerate}
\item When $n=0,1$ the set $\mathcal{M}^c_n(q,p)$ is either empty or a smooth manifold
      without boundary.
\item For $n > 1$ the set $\mathcal{M}^c_n(q,p)$ is either empty or a smooth manifold with
      corners.
\item The set $\mathcal{M}^c(q,p)$ is either empty or a smooth manifold without boundary.
\end{enumerate}
In each case the dimension of the manifold is $\lambda_q - \lambda_p - 1$.
When $M$ is orientable and $C_k$ is orientable for all $k=1,\ldots, l$, the above manifolds
are orientable.
\end{theorem}

At first glance, it may seem strange that $\mathcal{M}^c_n(q,p)$ is a smooth manifold
with corners for $n > 1$ whereas 
$$
\mathcal{M}^c(q,p) = \bigcup_{n \in \mathbb{Z}_+} \mathcal{M}^c_n(q,p)
$$
is a smooth manifold without boundary.  However, the proof of the above theorem in
\cite{BanCas} shows that the smooth manifolds with corners $\mathcal{M}^c_n(q,p)$ glue
together to form the manifold without boundary $\mathcal{M}^c(q,p)$, similar to the way 
that the manifold with boundary $[0,\infty)$ can be glued to $(-\infty,0)$ to create the
manifold without boundary $(-\infty,\infty)$.  


\subsection*{Compactness for moduli spaces of cascades}
In order to define a boundary operator by counting the number of cascades between two
critical points, the compactness properties of the manifold $\mathcal{M}^c(q,p)$
must be addressed. As one might expect, $\mathcal{M}^c(q,p)$ in general won't be compact
unless $\lambda_q - \lambda_p = 1$, because a sequence in $\mathcal{M}^c(q,p)$ may converge
to a broken flow line with cascades from $q$ to $p$. However, the precise definition
of a ``broken flow line with cascades'' turns out to be more subtle than the definition of
a ``broken gradient flow line''.

For a Morse-Bott function $f:M \rightarrow \mathbb{R}$, a broken gradient flow
line is simply a concatenation of gradient flow lines. As such, a broken gradient flow line
can be represented by an $n$-tuple $(x_1,\ldots ,x_n)$ where $x_k$ is a gradient
flow line of $f$ for all $1 \leq k \leq n$ and $\lim_{t \rightarrow \infty} x_k(t) =
\lim_{t \rightarrow -\infty} x_{k+1}(t)$ for all $1 \leq k \leq n-1$. The second
condition can be interpreted as saying that the time spent flowing along each
intermediate critical submanifold is $0$, and hence there is an obvious identification
of the broken gradient flow line represented by $(x_1,\ldots ,x_n)$ with the (non-broken)
flow line with $n$ cascades $((x_k)_{1 \leq k \leq n}, (t_k)_{1 \leq k \leq n-1})$ where 
$t_k = 0$ for all $1 \leq k \leq n-1$. 

This identification is compatible with the topology of the space of cascades $\mathcal{M}^c(q,p)$
and the topology of the space of broken gradient flow lines $\overline{\mathcal{M}}_f(C_j,C_i)$.
That is, suppose $\{\gamma_k\} \subset \mathcal{M}_f(C_j,C_i)$ is a sequence of unparameterized
gradient flow lines of a Morse-Bott-Smale function $f:M \rightarrow \mathbb{R}$, with 
$\partial_-(\gamma_k)\in W^u_{f_j}(q) \subseteq C_j$ and $\partial_+(\gamma_k) \in W^s_{f_i}(p)
\subset C_i$ for all $k$, that converges to a broken gradient flow line in 
$\overline{\mathcal{M}}_f(W^u_{f_j}(q), W^s_{f_i}(p))$ represented by $(x_1,\ldots ,x_n)$. Then
the proof of Theorem \ref{manifold} shows that the sequence $\{\gamma_k\}$, viewed as a subset of 
$\mathcal{M}^c(q,p)$, converges to the unparameterized flow line with cascades represented by 
$((x_k)_{1 \leq k \leq n}, (t_k)_{1 \leq k \leq n-1})$ where $t_k = 0$ for all $1 \leq k \leq n-1$.
Thus, the broken gradient flow lines from $W^u_{f_j}(q) \subseteq C_j$ to $W^s_{f_i}(p) \subseteq C_i$
are already included in the space of (unbroken) cascades $\mathcal{M}^c(q,p)$ in the sense that
the above identification induces an embedding $\overline{\mathcal{M}}_f(W^u_{f_j}(q), W^s_{f_i}(p))
\hookrightarrow \mathcal{M}^c(q,p)$ making the following diagram commute.
$$
\xymatrix{
\mathcal{M}_f(W^u_{f_j}(q),W^s_{f_i}(p)) \ar@{^{(}->}[dr] \ar@{^{(}->}[r] & 
  \overline{\mathcal{M}}_f(W^u_{f_j}(q),W^s_{f_i}(p)) \ar@{^{(}->}[d]\\
 & \mathcal{M}^c(q,p)
}
$$
Note: There are several equivalent ways of defining the topology on the space
of broken gradient flow lines $\overline{\mathcal{M}}_f(W^u_{f_j}(q), W^s_{f_i}(p))$, 
see Section 2 of \cite{HurFlo} or \cite{LetMor} for more details.

\smallskip
So, what is a ``broken flow line with cascades"? Upon further analysis it turns out that
a ``broken flow line with cascades'' should (roughly speaking) be a concatenation of
unparameterized flow lines with cascades that either flows along an intermediate critical
submanifold for infinite time or rests at an intermediate critical point of one of the Morse
functions on the critical submanifolds for infinite time. This description of the space
$\overline{\mathcal{M}}^c(q,p)$ of broken flow lines with cascades from $q$ to $p$ 
(and its topology) was made precise by Banyaga and Hurtubise in \cite{BanCas} by identifying
the set of broken flow lines with cascades with a set of compact subsets of a compact
metric space, whose topology is determined by the Hausdorff metric.

\begin{definition}
Let $(X,d)$ be a compact metric space and let $K_1$ and $K_2$ be nonempty closed subsets of $X$. 
The \textbf{Hausdorff distance} between $K_1$ and $K_2$ is defined to be
\begin{eqnarray*}
d_H(K_1,K_2) 
& = & \max \left\{\sup_{x_1 \in K_1} \inf_{x_2 \in K_2} d(x_1,x_2), \sup_{x_2 \in K_2}
      \inf_{x_1 \in K_1} d(x_1,x_2)\right\}\\
& = & \inf \left\{\varepsilon > 0 |\ K_1 \subseteq N_\varepsilon(K_2) \text{ and }
      K_2 \subseteq N_\varepsilon(K_1)   \right\}
\end{eqnarray*}
where $N_\varepsilon(K) = \bigcup_{y \in K} \{ x \in X |\ d(x,y) \leq \varepsilon \}$. 
\end{definition}

An unparameterized gradient flow line of a Morse-Bott function $f:M \rightarrow
\mathbb{R}$ can be identified with its image in $M$, and this image will be a compact
subset of $M$ diffeomorphic to $\overline{\mathbb{R}} = \mathbb{R} \cup \{\pm \infty\}$ as
long as we include the limits of the gradient flow in the image. However, an unparameterized
flow line with cascades may ``rest'' at an intermediate critical point, and hence the map
that sends an unparameterized flow line with cascades to its image might not be injective.
In order to get an injective map one needs to keep track of the times $t_k$ spent flowing
along or resting on the intermediate critical submanifolds. This leads to a continuous
injection
$$
\mathcal{M}^c(q,p) \hookrightarrow \mathcal{P}^c(M) \times \overline{\mathbb{R}}^l
$$
where $\mathcal{P}^c(M)$ denotes the space of all compact subsets of $M$ and $l$ is the number
of critical submanifolds.  All these ideas can then be extended to the space of broken flow lines
with cascades $\overline{\mathcal{M}}^c(q,p)$ by considering the images of broken gradient 
flow lines of the Morse functions $f_k:C_k \rightarrow \mathbb{R}$ on the critical submanifolds
and allowing the $t_k$ to be $\infty$. From this point of view, the topology on the space of
unparameterized broken flow lines with cascades is the topology the set inherits as a subspace
of $\mathcal{P}^c(M) \times \overline{\mathbb{R}}^l$, i.e. the topology determined by the
Hausdorff metric. 

\smallskip
In \cite{BanCas} Banyaga and Hurtubise used these ideas to prove the following theorem.

\begin{theorem}\label{compactification}
The space $\overline{\mathcal{M}}^c(q,p)$ of broken flow lines with cascades is compact, 
and there is an injection that restricts to a continuous embedding
$$
\mathcal{M}^c(q,p) \hookrightarrow \overline{\mathcal{M}}^c(q,p) \subset 
\mathcal{P}^c(M) \times \overline{\mathbb{R}}^l.
$$
Hence, every sequence of unparameterized flow lines with cascades from $q$ to $p$ has a 
subsequence that converges to a broken flow line with cascades from $q$ to $p$.
\end{theorem}

\noindent
The following fundamental property is a straightforward consequence of this theorem.

\begin{corollary}\label{finitenumber}
If $\lambda_q - \lambda_p = 1$, then $\mathcal{M}^c(q,p)$ is compact and hence a finite set. 
\end{corollary}


\subsection*{The cascade chain complex}
We are now in a position to use the moduli spaces $\mathcal{M}^c(q,p)$ to define a cascade
chain complex $(C_\ast^c(f),\partial_\ast^c)$ whose boundary operator is determined by
counting cascades. Let $Cr = \bigcup_{j=1}^l Cr(f_j)$ be the collection of critical points
of the Morse functions $f_j:C_j \rightarrow \mathbb{R}$, let $Cr_k \subseteq Cr$ be the
collection of critical points whose total index is $k$, and let $C_k^c(f)$ be the free
abelian group generated by the elements in $Cr_k$. We would like to define a boundary operator
$$
\partial^c_k:C_k^c(f) \rightarrow C_{k-1}^c(f)
$$
by counting the number of elements in $\mathcal{M}^c(q,p)$, where $q \in Cr_k$ and 
$p \in Cr_{k-1}$, either over $\mathbb{Z}_2$ or over $\mathbb{Z}$ with signs determined by
some orientations.

The approach taken in \cite{FraThe} is to count the cascades over $\mathbb{Z}_2$, which gives
a chain complex that computes the homology of $M$ with coefficients in $\mathbb{Z}_2$. One
reason for only counting the cascades mod $2$ in \cite{FraThe} is that the approach used 
there to construct the moduli spaces $\mathcal{M}^c(q,p)$ doesn't readily yield orientations
on the moduli spaces. In contrast, the approach used by Banyaga and Hurtubise to prove Theorem
\ref{manifold} shows that the moduli spaces are orientable when $M$ and the critical
submanifolds are orientable, and it is possible to define a coherent system of orientations
for the moduli spaces.

However, even though it would be possible to define a coherent system of orientations for
the moduli spaces $\mathcal{M}^c(q,p)$, the main theorem  in \cite{BanCas} is a 
correspondence theorem that says that when $\lambda_q - \lambda_p = 1$ there is a bijection
$$
\mathcal{M}^c(q,p) \leftrightarrow \mathcal{M}_{h_\varepsilon}(q,p)
$$
between the moduli space of cascades and the moduli space of gradient flow
lines  of the perturbed function 
$$
h_\varepsilon \stackrel{\mbox{def}}{=} f + \varepsilon \left( \sum_{j=1}^l \rho_j f_j \right)
$$
discussed in Section \ref{perturbsection} for $\varepsilon > 0$ sufficiently small. So, the
approach taken in \cite{BanCas} is to use the Correspondence Theorem to transfer the
orientations on $\mathcal{M}_{h_\varepsilon}(q,p)$ to $\mathcal{M}^c(q,p)$ and then define
the boundary operator $\partial^c_\ast$ over $\mathbb{Z}$ by counting cascades with signs
given by the induced orientations. 

This approach shows immediately that $\partial^c_\ast \circ \partial^c_\ast = 0$ and 
$$
H_\ast(C_\ast^c(f),\partial^c_\ast) \approx H_\ast(C_\ast(h_\varepsilon),
\partial_\ast^{h_\varepsilon}) \approx H_\ast(M;\mathbb{Z}).
$$
Moreover,  it proves that the chain complex defined using cascades is the same as the
Morse-Smale-Witten chain complex of the perturbed function $h_\varepsilon:M \rightarrow 
\mathbb{R}$, i.e. the generators of the two chain complexes are the same and the boundary
operators agree up to sign. This is a much stronger result than the statement that the
two chain complexes compute the same homology.

\smallskip\noindent
\textbf{Remark.}
It should be noted that Bourgeois and Oancea used a similar approach to orienting
moduli spaces of cascades in the context of symplectic homology \cite{BouAne} \cite{BouSym}.
That is, they proved a correspondence theorem between moduli spaces of cascades and
moduli spaces of flow lines of a perturbed function, and then they oriented the moduli 
spaces of cascades using their Correspondence Theorem. In particular, see Theorem 3.7
(Correspondence Theorem) in \cite{BouSym} and the discussion that follows.


\subsection*{Proving the Correspondence Theorem using the Exchange Lemma}
The proof of the Correspondence Theorem in \cite{BanCas} has several steps. Starting with
a Morse-Bott function $f:M \rightarrow \mathbb{R}$, a Riemannian metric $g$ on $M$ such that
$f$ satisfies the Morse-Bott-Smale transversality condition with respect to $g$, and
an $\varepsilon > 0$ small enough so that a list of conditions are met, Banyaga and Hurtubise
first show that there exists a small perturbation of the metric to a metric $\tilde{g}$
such that $h_{\varepsilon'}:M \rightarrow \mathbb{R}$ satisfies the Morse-Smale 
transversality condition with respect to $\tilde{g}$ for all $0 < \varepsilon' \leq 
\varepsilon$.  The perturbation can be chosen small enough so that $f$ satisfies the
Morse-Bott-Smale transversality condition with respect to $\tilde{g}$ and the hypotheses 
of Theorem \ref{manifold} still hold. Hence, there exists a metric $\tilde{g}$ such
that moduli spaces of cascades are defined and $h_{\varepsilon'}$ satisfies the 
Morse-Smale transversality condition for all $0 < \varepsilon' \leq \varepsilon$.

\smallskip\noindent
Banyaga and Hurtubise then prove the following lemma.

\begin{lemma}\label{constant}
Let $p,q \in \text{Cr}$ with $\lambda_q - \lambda_p = 1$, and let  $0 < \varepsilon' \leq
\varepsilon$.  If $h_{\varepsilon'}:M \rightarrow \mathbb{R}$ and $h_\varepsilon:M 
\rightarrow \mathbb{R}$ are Morse-Smale with respect to the same Riemannian metric, 
then the number of gradient flow lines of $h_{\varepsilon'}$ from $q$ to $p$ is equal 
to the number of gradient flow lines of $h_\varepsilon$ from $q$ to $p$.
\end{lemma}

\noindent
This lemma shows that with respect to the perturbed metric $\tilde{g}$ from above
there is a trivial cobordism
$$
\mathcal{M}_{h_\varepsilon}(q,p) \times (0,\varepsilon]
$$
such that
$$
\mathcal{M}_{h_\varepsilon}(q,p) \times \{\varepsilon'\} \approx \mathcal{M}_{h_{\varepsilon'}}
(q,p)
$$
for all $0 < \varepsilon' \leq \varepsilon$.  The next step is to analyze what
happens as $\varepsilon' \rightarrow 0$. This is sometimes referred to as
``degenerating the asymptotics''.

\begin{lemma}\label{degenerating}
Let $\{\varepsilon_\nu\}_{\nu=1}^\infty$ be a decreasing sequence such that $0 < \varepsilon_\nu
\leq \varepsilon$ for all $\nu$ and $\lim_{\nu \rightarrow \infty} \varepsilon_\nu = 0$. Let 
$q,p \in \text{Cr}$, and suppose that $\gamma_{\varepsilon_\nu} \in
\mathcal{M}_{h_{\varepsilon_\nu}}(q,p)$ for all $\nu$. Then there exists a broken flow line
with cascades $\gamma \in \overline{\mathcal{M}}^c(q,p)$ and a subsequence of
$\{Im(\gamma_{\varepsilon_\nu})\}_{\nu=1}^\infty$ that converges to $Im(\gamma)$ in the
Hausdorff topology.
\end{lemma}

\noindent
The proof of this lemma (which does not require $\lambda_q - \lambda_p = 1$) uses techniques
similar to those used to prove that the space of broken flow lines with cascades is compact
with respect to the Hausdorff topology (Theorem \ref{compactification}). If we use the
Hausdorff topology on
$$
\bigcup_{\varepsilon' \in (0,\varepsilon]} \mathcal{M}_{h_{\varepsilon'}}(q,p)
\subset \mathcal{P}^c(M)
$$
(where an element of $\mathcal{M}_{h_{\varepsilon'}}(q,p)$ is identified with its
image, including $q$ and $p$), then Lemma \ref{degenerating} says that the boundary of 
this space is contained in the union of the images of the cascades in $\mathcal{M}^c(q,p)$.
However, it is still possible that two distinct sequences of gradient flow lines from $q$
to $p$ might converge the same cascade as $\varepsilon' \rightarrow 0$ or some of the
cascades in $\mathcal{M}^c(q,p)$ might not be near any of the gradient flow lines of 
$h_{\varepsilon'}$, even when $\varepsilon' > 0$ is extremely small. So, the proof of
the Correspondence Theorem requires a much more detailed analysis than is provided by
Lemma \ref{degenerating}.

\smallskip
To conduct this more detailed analysis Banyaga and Hurtubise turned to the Exchange Lemma,
which can be viewed as a generalization of Palis' $\lambda$-Lemma. The $\lambda$-Lemma
applies to a critical point $p$ of a Morse-Smale system, and it says (roughly speaking) that
if $N$ is an invariant submanifold that intersects $W^s(p)$ transversally then $W^u(p)$ must
contain points that are close to $N$. The $\lambda$-Lemma is an essential tool for the dynamical
systems approach to studying compactified moduli spaces of Morse-Smale flows (cf. Sections 6.2
and 6.3 of \cite{BanLec}), and the Exchange Lemma allows the dynamical systems approach to be
extended to Morse-Bott-Smale systems.

\smallskip
The Exchange Lemma comes out of geometric singular perturbation theory, and it applies
to ``fast-slow'' dynamical systems.  Following the notation in \cite{JonGeo}, a fast-slow
system of differential equations in local coordinates is of the form
\begin{eqnarray*}
x' & = & f(x,y,\epsilon)\\
y' & = & \epsilon g(x,y,\epsilon)
\end{eqnarray*}
where $' = \frac{d}{dt}$, $x \in \mathbb{R}^n$, $y \in \mathbb{R}^l$, $\epsilon$ 
is a real parameter, and both $f$ and $g$ are $C^\infty$ (and hence bounded) on some
neighborhood of $0$. The $x$ coordinates are called the fast variables and the
$y$ coordinates are called the slow variables because in the limit as $\epsilon 
\rightarrow 0$ we have
\begin{eqnarray*}
x' & = & f(x,y,0)\\
y' & = & 0
\end{eqnarray*}
where the $x$ coordinates can vary but $y$ remains constant. Alternately, when 
$\epsilon \neq 0$ is close to $0$, $y'$ is close to $0$ and the $y$ coordinates change
slowly, whereas the $x$ coordinates can change more quickly.

In the setup contained in \cite{BanCas}, each critical submanifold has a neighborhood with
coordinates $(u,v,w)$ coming from the Morse-Bott Lemma, where the $u$ coordinates
are the coordinates along the critical submanifold and the $(v,w)$ coordinates are the
coordinates in the directions normal to the critical submanifold. The Morse function on
the critical submanifold depends only on the $u$ coordinates, which are the slow variables,
and the Morse-Bott function depends only on the $(v,w)$ coordinates, which are the fast
variables. In fact, the Riemannian metric is chosen so that on a neighborhood of the
critical submanifold
$$
\nabla h_\varepsilon = \nabla f + \varepsilon \nabla f_j
$$
where $\nabla f \perp \nabla f_j$.  Thus, the gradient flow equation of
the Morse-Smale function $\nabla h_\varepsilon$ in the local coordinates $(u,v,w)$
near the critical submanifold is
\begin{eqnarray*}
(v',w') & = & (\nabla f)(v,w)\\
u'      & = & \varepsilon (\nabla f_j)(u)
\end{eqnarray*}
which is a fast-slow system.

Several versions of the Exchange Lemma with various levels of generality have been proved
by many different authors, cf. \cite{JonGen} \cite{LiuGeo} \cite{SchExcI} \cite{SchExcII}.
The lemma gives a relationship between the dynamics of a fast-slow system when
$\varepsilon \neq  0$ and the dynamics of the system when $\varepsilon = 0$.
Roughly speaking, the lemma says that a manifold $M_0$ that is transverse to the stable
manifold $W^s_0(C)$ of a normally hyperbolic critical submanifold $C$ of the system with
$\epsilon = 0$ will have points that flow forward in time under the fast-slow
system with $\varepsilon\neq 0$ to be near subsets of the unstable manifold $W^u_0(C)$
of the system with $\varepsilon = 0$.

With respect to the setup in \cite{BanCas}, we have a flow line with $n$ cascades
$$
\left((x_k)_{1\leq k\leq n},(t_k)_{1 \leq k \leq n-1}  \right)
$$
with intermediate critical submanifolds $C_{j_1}, \ldots , C_{j_n}$ and local coordinates
$(u,v,w)$ near an intermediate critical submanifold $C_{j_k}$.  In the local coordinates
the critical submanifold $C_{j_k}$ consists of the $u$ components $\{(u,0,0)\}$,
the stable manifold $W^s_f(C_{j_k})$ is given by the $(u,w)$ components $\{(u,0,w)\}$,
and the unstable manifold $W^s_f(C_{j_k})$ is given by the $(u,v)$ components $\{(u,v,0)\}$.
Away from the critical submanifolds the gradient flow lines $x_k(t)$ of $f = h_0$ and 
$x_k^\varepsilon(t)$ of $h_\varepsilon$ agree. However, near the critical submanifold
$C_{j_k}$ the gradient of $h_\varepsilon$ with $\varepsilon \neq 0$ may be nonzero in the
$u$ components, whereas the gradient of $f$ is zero in the $u$ components. So, near the
critical submanifold the gradient flow line $x_k^\varepsilon(t)$ of $h_\varepsilon$ can 
diverge from the gradient flow line $x_k(t)$ of $f=h_0$.

\begin{figure}[h]
\includegraphics{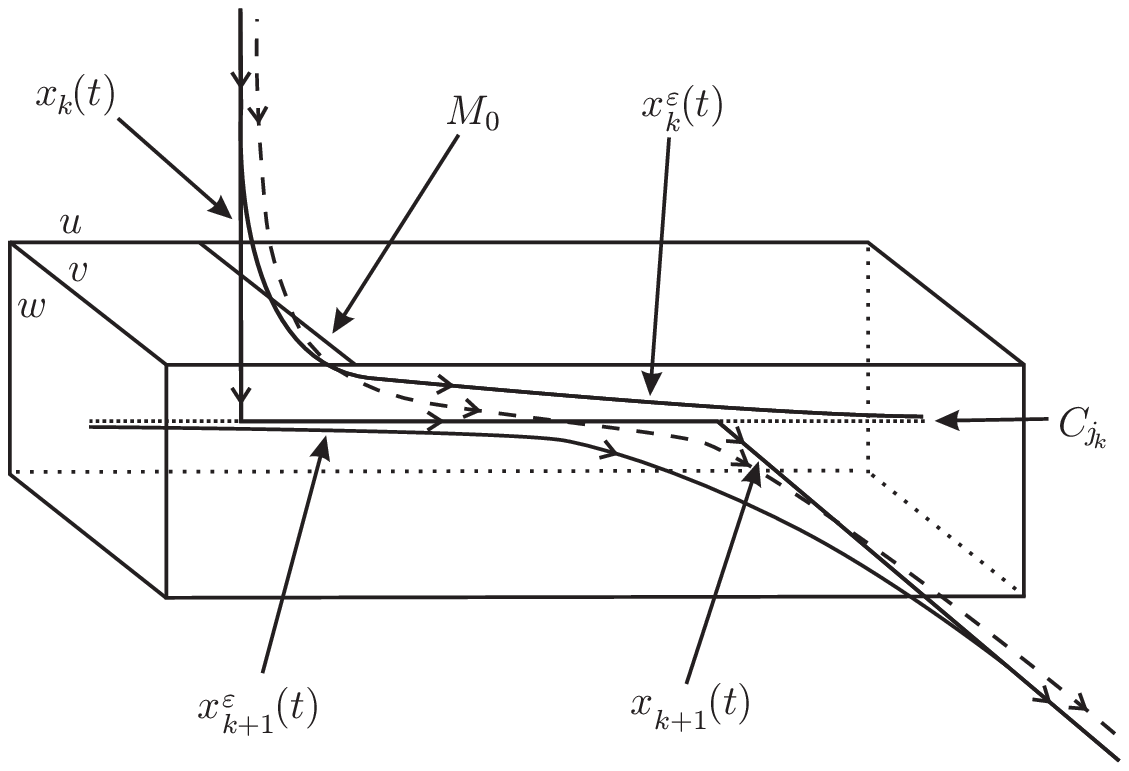}
\end{figure}

The Exchange Lemma says that if $M_0$ intersects $W^s_f(C_{j_k})$ transversally near
the image of $x_k(t)$, then there must be points in $M_0$ that flow forward in time under
the gradient flow of $h_\varepsilon$ with $\varepsilon \neq 0$ to be near the image of
$x_{k+1}(t)$. Thus, there is a gradient flow line of $h_\varepsilon$ passing through
$M_0$ whose image is near the image of the cascade $(x_k,x_{k+1},t_k)$. This is shown in the
diagram where the gradient flow line of $h_\varepsilon$ near the image of the cascade is
the dashed curve lying above the unstable manifold $\{(u,v,0)\}$.

\smallskip
Using these ideas, Banyaga and Hurtubise proved the main theorem in \cite{BanCas},
which implies that the cascade chain complex is the same as the Morse-Smale-Witten
chain complex of $h_\varepsilon$ up to sign.

\begin{theorem}[Correspondence of Moduli Spaces]
Let $p,q \in \text{Cr}(h_\varepsilon)$ with $\lambda_q - \lambda_p = 1$.
For any sufficiently small $\varepsilon>0$ there is a bijection between unparameterized
cascades and unparameterized gradient flow lines of the Morse-Smale function
$h_\varepsilon:M \rightarrow \mathbb{R}$ between $q$ and $p$,
$$
\mathcal{M}^c(q,p) \leftrightarrow \mathcal{M}_{h_\varepsilon}(q,p).
$$
\end{theorem}

\noindent
The Correspondence Theorem allows us to identify the space of cascades
$\mathcal{M}^c(q,p)$ with the left side boundary of the trivial cobordism
$$
\mathcal{M}_{h_\varepsilon}(q,p) \times [0,\varepsilon],
$$
which will have the opposite orientation as the right side boundary.

\begin{corollary}[Correspondence of Chain Complexes]
For $\varepsilon >0$ sufficiently small, the Morse-Smale-Witten chain complex
$(C_\ast(h_\varepsilon),\partial_\ast^{h_\varepsilon})$ associated to the perturbation
$$
h_\varepsilon = f + \varepsilon \left( \sum_{j=1}^l \rho_j f_j \right)
$$
of a Morse-Bott function $f:M \rightarrow \mathbb{R}$ is the same as the cascade
chain complex $(C_\ast^c(f),\partial^c_\ast)$ up to sign. That is, the chain groups of both complexes
have the same generators and $\partial^c_\ast = - \partial_\ast^{h_\varepsilon}$.
\end{corollary}


\section{The Morse-Bott multicomplex}\label{multicomplexsection}

The approaches discussed in the previous sections require choosing auxiliary Morse
functions on the critical submanifolds in order to define a chain complex generated by
the critical points of the chosen Morse functions. The approach discussed in this
section does not involve choosing any auxiliary Morse functions. Instead, the chain groups
are generated by singular topological chains on the critical submanifolds. Keeping track
of the degrees of the singular topological chains, the Morse-Bott indexes of the critical
submanifolds, and homomorphisms defined using moduli spaces of gradient flow lines between
the critical submanifolds leads to an algebraic structure known as a multicomplex, which
generalizes the notion of a double complex.


\subsection*{Multicomplexes and assembled chain complexes}

\begin{definition}\label{multicomplex}
Let $R$ be a principal ideal domain. A first quadrant \textbf{multicomplex} $X$ is 
a bigraded $R$-module $\{X_{p,q}\}_{p,q \in \mathbb{Z}_+}$ with
differentials
$$
\textmc{d}_j:X_{p,q} \rightarrow X_{p-j,q+j-1} \quad \text{ for all } j = 0,1,\ldots
$$
that satisfy
$$
\sum_{i+j = n} \textmc{d}_i \textmc{d}_j = 0 \quad \text{ for all }n.
$$
A first quadrant multicomplex such that $\textmc{d}_j = 0$ for all $j \geq 2$ is called
a \textbf{double complex} (or a \textbf{bicomplex}).
\end{definition}

\noindent
A first quadrant multicomplex looks similar to a spectral sequence, but the differentials
are all defined on the $0^\text{th}$ page and we may have $\textmc{d}_j \circ \textmc{d}_j
\neq 0$ when $j > 0$.

\begin{figure}[h]
A First Quadrant Multicomplex
$$
\xymatrix{
\vdots & \vdots & \vdots & \vdots \\
X_{0,3} \ar[d]^(.4){\textmc{d}_0} & X_{1,3} \ar[d]^(.4){\textmc{d}_0} \ar[l]_{\textmc{d}_1} & X_{2,3} \ar[d]^(.4){\textmc{d}_0} \ar[l]_{\textmc{d}_1} & X_{3,3} \ar[d]^(.4){\textmc{d}_0}
\ar[l]_{\textmc{d}_1} & \cdots\\
X_{0,2} \ar[d]^(.4){\textmc{d}_0} & X_{1,2} \ar[d]^(.4){\textmc{d}_0} \ar[l]_{\textmc{d}_1} & X_{2,2} \ar[d]^(.4){\textmc{d}_0} \ar[l]_{\textmc{d}_1} \ar[llu]|(.37){\textmc{d}_2} |!{[ul];[l]}\hole & 
X_{3,2} \ar[d]^(.4){\textmc{d}_0} \ar[l]_{\textmc{d}_1} 
\ar[llu]|(.37){\textmc{d}_2} |!{[ul];[l]}\hole & \cdots\\
X_{0,1} \ar[d]^(.4){\textmc{d}_0} & X_{1,1} \ar[d]^(.4){\textmc{d}_0} \ar[l]_{\textmc{d}_1} & X_{2,1} \ar[d]^(.4){\textmc{d}_0} \ar[l]_{\textmc{d}_1} \ar[llu]|(.37){\textmc{d}_2} |!{[ul];[l]}\hole & 
X_{3,1} \ar[d]^(.4){\textmc{d}_0} \ar[l]_{\textmc{d}_1} \ar[llu]|(.37){\textmc{d}_2} |!{[ul];[l]}\hole 
  \ar@{.>}[llluu]_>(.2){\textmc{d}_3} & \cdots \\
X_{0,0}  & X_{1,0} \ar[l]_{\textmc{d}_1} & X_{2,0} \ar[l]_{\textmc{d}_1} \ar[llu]|(.37){\textmc{d}_2} |!{[ul];[l]}\hole & X_{3,0} \ar[l]_{\textmc{d}_1} \ar[llu]|(.37){\textmc{d}_2} |!{[ul];[l]}\hole \ar@{.>}[llluu]_>(.2){\textmc{d}_3} & \cdots
}
$$
\end{figure}

A multicomplex can be \textbf{assembled} to form a filtered chain complex 
$((CX)_\ast,\partial_\ast)$ by summing along the diagonals. That is, if we define
$$
(CX)_k \equiv \bigoplus_{p+q = k} X_{p,q}
$$
and $\partial_k = \textmc{d}_0 \oplus \cdots \oplus \textmc{d}_k$ for all $k \in \mathbb{Z}_+$,
then the relations in Definition \ref{multicomplex} imply that $\partial_k \circ \partial_{k+1} = 0$.

\smallskip\noindent
Note: The chain complex $((CX)_\ast,\partial_\ast)$ has a filtration given by
$$
F_s(CX)_k\ \equiv \bigoplus_{\stackrel{\scriptstyle p+q = k}{p\leq s}} X_{p,q}
$$
which determines a spectral sequence.  However, the differentials in this spectral sequence
are not necessarily induced from the differentials $\textmc{d}_j$ when $j \geq 2$ \cite{HurMul}.
\nocite{RuaBot}

\begin{figure}[h]
The Assembled Chain Complex
$$
\xymatrix{
\ddots & \vdots \\
\cdots & X_{3,0} \ar[r]^(.53){\textmc{d}_0}  \ar[dr]^{\textmc{d}_1} 
\ar[ddr]|(.33){\textmc{d}_2} |!{[d];[dr]}\hole \ar@{.>}[dddr]|(.25){\textmc{d}_3} & 0 \\
\cdots & X_{2,1} \ar@{}[u]|{\oplus} \ar[r]^(.53){\textmc{d}_0} \ar[dr]^{\textmc{d}_1} 
\ar[ddr]|(.33){\textmc{d}_2} |!{[d];[dr]}\hole & X_{2,0} \ar@{}[u]|{\oplus} \ar[r]^(.53){\textmc{d}_0} 
\ar[dr]^{\textmc{d}_1} \ar[ddr]|(.33){\textmc{d}_2} |!{[d];[dr]}\hole & 0 & & \\
\cdots & X_{1,2} \ar@{}[u]|{\oplus} \ar[r]^(.53){\textmc{d}_0} \ar[dr]^{\textmc{d}_1} & X_{1,1} 
\ar@{}[u]|{\oplus} \ar[r]^(.53){\textmc{d}_0} \ar[dr]^{\textmc{d}_1} & X_{1,0} \ar@{}[u]|{\oplus} 
\ar[r]^{\textmc{d}_0} \ar[dr]^{\textmc{d}_1} & 0 & \\
\cdots & X_{0,3} \ar@{}[u]|{\oplus} \ar[r]^(.53){\textmc{d}_0} & X_{0,2} \ar@{}[u]|{\oplus} 
\ar[r]^(.53){\textmc{d}_0} & X_{0,1} \ar@{}[u]|{\oplus} \ar[r]^(.53){\textmc{d}_0} & X_{0,0} \ar@{}[u]|{\oplus}
\ar[r]^-{\textmc{d}_0} & 0\\
\cdots &  {(CX)}_3 \ar@{}[u]|{\|} \ar[r]^{\partial_3} & {(CX)}_2 \ar@{}[u]|{\|} 
\ar[r]^{\partial_2} & {(CX)}_1 \ar@{}[u]|{\|} \ar[r]^{\partial_1} & {(CX)}_0 \ar@{}[u]|{\|}
\ar[r]^-{\partial_0} & 0 \ar@{}[u]|{\|}
}
$$
\end{figure}


\subsection*{A heuristic view of the Morse-Bott multicomplex}
Let $f:M \rightarrow \mathbb{R}$ be a Morse-Bott-Smale function on an $m$-dimensional
compact smooth closed Riemannian manifold $M$, and let $B_i \subseteq Cr(f)$ be
the union of the critical submanifolds of Morse-Bott index $i$ for $i=0,\ldots ,m$.
The compactified moduli space $\overline{\mathcal{M}}(B_i,B_{i-j})$ of broken gradient flow
lines of $f$ from $B_i$ to $B_{i-j}$ for $j=1,\ldots, i$ is a smooth manifold with corners
and the beginning point map
$$
\partial_-: \overline{\mathcal{M}}(B_i,B_{i-j}) \rightarrow B_i
$$
is a submersion and a stratum submersion (cf. Corollary 5.20 of \cite{BanMor}). Thus, every
smooth map $\sigma:P \rightarrow B_i$ from a smooth manifold with corners $P$ is transverse
and stratum transverse to $\partial_-$, and the fibered product
$P \times_{B_i} \overline{\mathcal{M}}(B_i,B_{i-j})$ of $\sigma$ and $\partial_-$
over $B_i$ is a smooth manifold with corners (cf. Lemma 5.21 of \cite{BanMor}).
$$
\xymatrix{
P \times_{B_i} \overline{\mathcal{M}}(B_i,B_{i-j}) \ar@{-->}[r] \ar@{-->}[d] & 
  \overline{\mathcal{M}}(B_i,B_{i-j}) \ar[d]^{\partial_-}\\
P \ar[r]^{\sigma} & B_i
}
$$
(Similar spaces were used in the proof of Theorem \ref{manifold} on moduli spaces of cascades).
Composing the projection map $\pi_2$ onto the second component of 
$P \times_{B_i} \overline{\mathcal{M}}(B_i,B_{i-j})$ with the endpoint map 
$\partial_+:\overline{\mathcal{M}}(B_i,B_{i-j}) \rightarrow B_{i-j}$ gives a map
$$
P \times_{B_i} \overline{\mathcal{M}}(B_i,B_{i-j}) 
\stackrel{\pi_2}{\longrightarrow} \overline{\mathcal{M}}(B_i,B_{i-j})
\stackrel{\partial_+}{\longrightarrow} B_{i-j}.
$$
Moreover, if $P$ has dimension $p$, then $P \times_{B_i} \overline{\mathcal{M}}(B_i,B_{i-j})$
has dimension $p+j-1$, which is independent of the dimension of the connected components in
$B_i$ and $B_{i-j}$. 

Up to this point, the discussion has been rigorous.  We will now make explicit an
unwarranted assumption that has been assumed implicitly by other authors (cf. \cite{FukFlo}
\cite{LiuOnt}). If the above fibered product had a preferred finite triangulation, then
summing over the restrictions of the above map to the simplices making up the finite
triangulation would define a singular chain $\partial_j(\sigma)$ in $B_{i-j}$. Moreover,
if every smooth manifold with corners under consideration came with a preferred finite
triangulation (or cubulation), then this fibered product construction would define a
homomorphism $\partial_j:S_p(B_i) \rightarrow S_{p+j-1}(B_{i-j})$ from the singular 
$p$-chains on $B_i$ to the singular $p+j-1$-chains on $B_{i-j}$ (or singular cubical chains 
if we were given preferred finite cubulations). These maps would then yield the following,
where $\partial_0$ is comes from the usual singular boundary operator.

\begin{figure}[h]
Heuristic View of the Morse-Bott Multicomplex
$$
\xymatrix{
\vdots & \vdots & \vdots & \vdots \\
S_3(B_0) \ar[d]^(.4){\partial_0} & S_3(B_1) \ar[d]^(.4){\partial_0} \ar[l]_{\partial_1} & 
S_3(B_2) \ar[d]^(.4){\partial_0} \ar[l]_{\partial_1} & S_3(B_3) \ar[d]^(.4){\partial_0}
\ar[l]_{\partial_1} & \cdots\\
S_2(B_0) \ar[d]^(.4){\partial_0} & S_2(B_1) \ar[d]^(.4){\partial_0} \ar[l]_{\partial_1} & 
S_2(B_2) \ar[d]^(.4){\partial_0} \ar[l]_{\partial_1} \ar[llu]|(.37){\partial_2} |!{[ul];[l]}\hole & S_2(B_3) \ar[d]^(.4){\partial_0} \ar[l]_{\partial_1} 
\ar[llu]|(.37){\partial_2} |!{[ul];[l]}\hole & \cdots\\
S_1(B_0) \ar[d]^(.4){\partial_0} & S_1(B_1) \ar[d]^(.4){\partial_0} \ar[l]_{\partial_1} & 
S_1(B_2) \ar[d]^(.4){\partial_0} \ar[l]_{\partial_1} \ar[llu]|(.37){\partial_2} |!{[ul];[l]}\hole & S_1(B_3) \ar[d]^(.4){\partial_0} \ar[l]_{\partial_1} \ar[llu]|(.37){\partial_2}
 |!{[ul];[l]}\hole \ar@{.>}[llluu]_>(.2){\partial_3} & \cdots \\
S_0(B_0) & S_0(B_1) \ar[l]_{\partial_1} & S_0(B_2) \ar[l]_{\partial_1} \ar[llu]|(.37){\partial_2} |!{[ul];[l]}\hole & S_0(B_3) \ar[l]_{\partial_1} \ar[llu]|(.37){\partial_2} |!{[ul];[l]}\hole \ar@{.>}[llluu]_>(.2){\partial_3} & \cdots
}
$$
\end{figure}

Of course, smooth manifolds with corners don't usually come with preferred triangulations,
and there is no preferred (or induced) finite triangulation on the fibered product of
finitely triangulated spaces (cf. Example 5.17 of \cite{BanMor}). Still, it might
be possible to pick finite triangulations on all the (uncountably many) spaces under
consideration, prove that the relations in Definition \ref{multicomplex} hold with respect
to the chosen triangulations, and then show that the homology of the resulting assembled
chain complex is independent of the chosen triangulations.  However, there seem to be many
technical difficulties involved with making this approach rigorous on the level of chains.
Fortunately, by expanding the collection of allowed domains for the singular chains
it is possible to construct a Morse-Bott multicomplex without choosing any
triangulations.


\subsection*{The Banyaga-Hurtubise approach to the Morse-Bott multicomplex}

Singular homology is usually defined using maps from the standard $k$-simplex $\Delta^k$.
However, other equivalent versions of singular homology have been defined using maps from
domains other than $\Delta^k$.  For instance, there is singular cubical homology, which 
is based on maps from the unit $k$-cube $I^k$ \cite{MasABa}, and there is also a version
of singular homology based on maps from permutahedra \cite{UmbDia}. In order to create a
singular homology theory that allows for even more general domains Banyaga and Hurtubise 
make the following definitions in Section 4 of \cite{BanMor}.

\smallskip
For each integer $p \geq 0$ fix a set $C_p$ of topological spaces, and
let $S_p$ be the free abelian group generated by the elements of $C_p$, 
i.e. $S_p = \mathbb{Z}[C_p]$. Set $S_p = \{0\}$ if $p<0$ or $C_p = \emptyset$. 

\begin{definition}\label{topologicalchainsdef}
A \textbf{boundary operator} on the collection $S_\ast$ of groups $\{S_p\}$ 
is a homomorphism $\partial_p:S_p \rightarrow S_{p-1}$ such that 
\begin{enumerate}
\item For $p \geq 1$ and $P \in C_p\subseteq S_p$, $\partial_p(P) = \sum_k n_k P_k$
      where $n_k = \pm 1$ and $P_k\in C_{p-1}$ is a subspace of $P$ for all $k$.
\item $\partial_{p-1}\circ \partial_p:S_p \rightarrow S_{p-2}$ is zero.
\end{enumerate}
The pair $(S_\ast,\partial_\ast)$ is called a \textbf{chain complex of abstract 
topological chains}, and elements of $S_p$ are called \textbf{abstract topological chains} 
of \textbf{degree} $p$.
\end{definition}

\begin{definition}\label{singulartopologicaldef}
Let $B$ be a topological space and $p \in \mathbb{Z}_+$. A
\textbf{singular $C_p$-space} in $B$ is a continuous map $\sigma:P 
\rightarrow B$ where $P\in C_p$, and the \textbf{singular $C_p$-chain group} 
$S_p(B)$ is the free abelian group generated by the singular $C_p$-spaces. 
Define $S_p(B) = \{0\}$ if $S_p = \{0\}$ or $B =\emptyset$.
Elements of $S_p(B)$ are called \textbf{singular topological chains}
of \textbf{degree} $p$. 
\end{definition}

\noindent
For $p \geq 1$ there is a boundary operator $\partial_p:S_p(B) \rightarrow S_{p-1}(B)$
induced from the boundary operator $\partial_p:S_p \rightarrow S_{p-1}$. If 
$\sigma:P \rightarrow B$ is a singular $C_p$-space in $B$, then $\partial_p(\sigma)$
is given by the formula
$$
\partial_p(\sigma) = \sum_k n_k \sigma|_{P_k}
$$ 
where 
$$
\partial_p(P) = \sum_k n_k P_k.
$$
The pair $(S_\ast(B),\partial_\ast)$ is called a \textbf{chain complex of
singular topological chains}.

\smallskip\noindent
{\bf Example: Singular $N$-cube chains.}
Pick some large positive integer $N$ and let $I^N = \{(x_1,\ldots ,x_N) \in
\mathbb{R}^N |\ 0\leq x_j \leq 1, \ j=1,\ldots ,N \}$ denote the unit
$N$-cube.  For every $0 \leq p \leq N$ let $C_p$ be the set consisting of
the faces of $I^N$ of dimension $p$, i.e. subsets of $I^N$ where $p$ of the
coordinates are free and the rest of the coordinates are fixed to be either
$0$ or $1$. For every $0 \leq p \leq N$ let $S_p$ be the free abelian group
generated by the elements of $C_p$.
For $P \in C_p$  define
$$
\partial_p(P) = \sum_{j=1}^p (-1)^j \left[ P|_{x_j = 1} - P|_{x_j = 0} \right]
\in S_{p-1}
$$
where $x_j$ denotes the $j^{\mbox{th}}$ free coordinate of $P$. 
It is easy to show that $\partial_{p-1} \circ \partial_p = 0$, and hence the
faces of $I^N$ are abstract topological chains. Thus, a continuous map $\sigma_P:P
\rightarrow B$ from a face of $I^N$ of dimension $p$ into a topological space $B$
is a singular $C_p$-space in $B$, and the boundary operator applied to $\sigma_P$ is
$$
\partial_p(\sigma_P) = \sum_{j=1}^p (-1)^j \left[\sigma_P|_{x_j = 1} -
\sigma_P|_{x_j = 0}\right] \in S_{p-1}(B)
$$
where $\sigma_P|_{x_j = 0}$ denotes the restriction $\sigma_P: P|_{x_j=0}
\rightarrow B$ and $\sigma_P|_{x_j = 1}$ denotes the restriction
$\sigma_P:P|_{x_j = 1} \rightarrow B$. 

For instance, if $p=N=2$ the abstract topological chain $I^2$ has boundary,

\begin{center}
\includegraphics{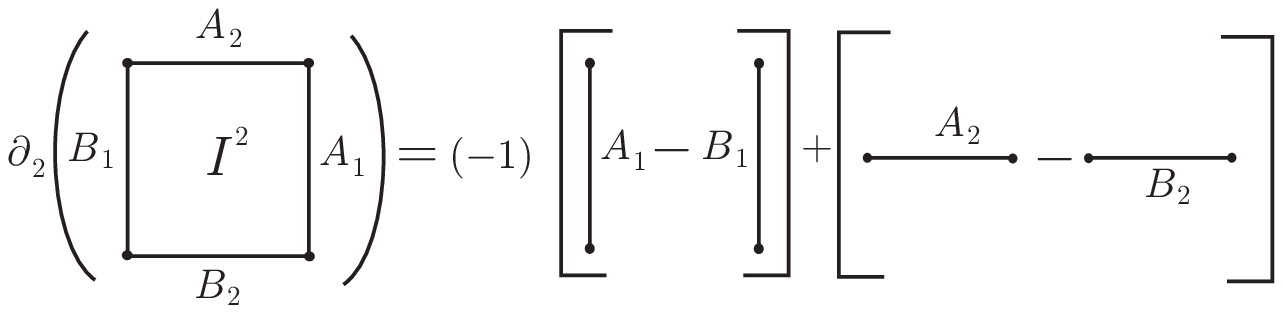}
\end{center}

\noindent
and the singular $C_2$-space $\sigma:I^2 \rightarrow B$ has boundary 
$$
\partial_2(\sigma) = (-1)[\sigma|_{A_1} - \sigma|_{B_1}] + [\sigma|_{A_2} - \sigma|_{B_2}].
$$
Note that this differs from the usual boundary operator on singular cubical chains
because there are several different domains of the same dimension. Normally,
singular homology is defined by picking a unique domain in each dimension and then
defining the boundary operator using inclusion maps. For instance, the boundary operator 
on singular cubes found in \cite{MasABa} is defined using the following inclusion maps 
when $p=2$.

\begin{center}
\includegraphics{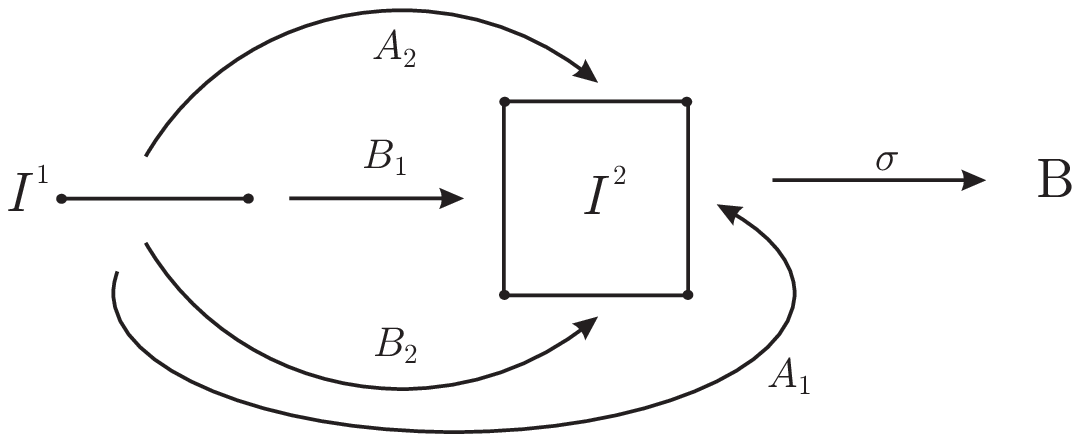}
\end{center}

To account for the multiple domains in each dimension, Banyaga and Hurtubise define
degeneracy relations in the form of a subgroup $D_p(B) \subseteq S_p(B)$ that identifies
maps that are ``essentially'' the same.  They then prove the following theorem.

\begin{theorem}[Singular $N$-Cube Chain Theorem]\label{cubehomology}
The boundary operator for singular $N$-cube chains 
$\partial_p:S_p(B) \rightarrow S_{p-1}(B)$ descends to a homomorphism 
$$
\partial_p:S_p(B)/D_p(B) \rightarrow S_{p-1}(B)/D_{p-1}(B),
$$
and
$$
H_p(S_\ast(B)/D_\ast(B), \partial_\ast) \approx H_p(B;\mathbb{Z})
$$ 
for all $p < N$.
\end{theorem}


\subsection*{Fibered products and moduli spaces as abstract topological chains}
Since most of the homomorphisms in the Morse-Bott multicomplex are defined
using fibered products of compactified moduli space of gradient flow lines,
the next step is to show that the compactified moduli spaces of gradient
flow lines of a Morse-Bott-Smale function are abstract topological chains and
a boundary operator on abstract topological chains extends to fibered products.
\smallskip

Let $f:M \rightarrow \mathbb{R}$ be a Morse-Bott-Smale function on an $m$-dimensional
compact smooth closed Riemannian manifold $M$, and let $B_i \subseteq Cr(f)$ be
the union of the critical submanifolds of Morse-Bott index $i$ for $i=0,\ldots ,m$. To 
simplify the notation in the following we will drop the subscript on $\partial$
and assume that for each $i=0,\ldots ,m$ the components of $B_i$ are all of the same
dimension. In general one needs to group the components by their dimension and then define
the degree and boundary operator on each group.

\begin{definition}\label{moduliboundary}
Let $B_i$ be the set of critical points of index $i$. For any $j=1,\ldots ,i$
the degree of $\overline{\mathcal{M}}(B_i,B_{i-j})$ is defined to be $j+b_i -1$
and the boundary operator is defined to be
$$
\partial \overline{\mathcal{M}}(B_i,B_{i-j}) = (-1)^{i+b_i} \sum_{i-j<n<i}
\overline{\mathcal{M}}(B_i,B_n) \times_{B_n} \overline{\mathcal{M}}(B_n,B_{i-j})
$$
where $b_i = \mbox{dim }B_i$ and the fibered product is taken over the beginning
and endpoint maps $\partial_-$ and $\partial_+$. If $B_n = \emptyset$, then
$\overline{\mathcal{M}}(B_i,B_n) = \overline{\mathcal{M}}(B_n,B_{i-j}) = 0$.
\end{definition}

\noindent
In order to check that $\partial \circ \partial = 0$ we need to know how
$\partial$ extends to fibered products.

\begin{definition}\label{fiberedboundary}
Suppose that $\{C_p\}_{p \geq 0}$ is a collection of topological spaces
that is closed under the fibered product construction with respect
to some collection of maps, and assume that $(S_\ast, \partial_\ast)$ is
a chain complex of abstract topological chains based on some subset of
the collection $\{C_p\}_{p \geq 0}$.
If $\sigma_i = \sum_k n_{i,k} \sigma_{i,k} \in S_{p_i}(B)$ is defined for $i=1,2$ 
where $\sigma_{i,k}:P_{i,k} \rightarrow B$ is a singular $C_{p_i}$-space 
for all $k$, then the \textbf{fibered product} of $\sigma_1$ and 
$\sigma_2$ over $B$ is defined to be
$$
P_1 \times_B P_2 = \sum_{k,j} n_{1,k} n_{2,j}\ P_{1,k} \times_B P_{2,j}
$$
where $P_1 = \sum_k n_{1,k}P_{1,k} \in S_{p_1}$ and 
$P_2 = \sum_j n_{2,j}P_{2,j} \in S_{p_2}$. The \textbf{boundary operator} 
applied to the fibered product is defined to be
$$
\partial (P_1 \times_B P_2) = \partial P_1 \times_B P_2
+ (-1)^{p_1+b} P_1 \times_B \partial P_2.
$$
If $\sigma_i = 0$ for either $i=1$ or $2$, then we define $P_1 \times_B P_2 = 0$.
\end{definition}

\noindent
The following lemmas from Section 4 of \cite{BanMor} show that the fibered
product of abstract topological chains is an abstract topological chain and
the compactified moduli spaces of gradient flow lines of a Morse-Bott-Smale
function are abstract topological chains.  The signs $(-1)^{i+b_i}$ and $(-1)^{p_1+b}$
in Definitions \ref{moduliboundary} and \ref{fiberedboundary} are essential to
the proofs of these two lemmas.

\begin{lemma}\label{fiberedcomplex}
The fibered product of two singular topological chains is an abstract
topological chain, i.e. the boundary operator on fibered products is of
degree -1 and satisfies $\partial \circ \partial = 0$. Moreover, 
the boundary operator on fibered products is associative, i.e.  
$$
\partial((P_1 \times_{B_1} P_2) \times_{B_2} P_3) = 
\partial(P_1 \times_{B_1} (P_2 \times_{B_2} P_3)).
$$
\end{lemma}

\begin{lemma}\label{modulicomplex}
The degree and boundary operator for $\overline{\mathcal{M}}(B_i,B_{i-j})$
satisfy the axioms for abstract topological chains, i.e. the boundary operator
on compactified moduli spaces of gradient flow lines from Definition
\ref{moduliboundary} is of degree $-1$ and it satisfies $\partial \circ \partial = 0$.
\end{lemma}


\subsection*{The Morse-Bott-Smale multicomplex}
Fix some $N > \text{dim } M$, and for any $p\geq 0$ let $C_p$ be the set consisting of
the faces of $I^N$ of dimension $p$ and the connected components of degree $p$ of
fibered products of the form
$$
Q \times_{B_{i_1}} \overline{\mathcal{M}}(B_{i_1},B_{i_2}) \times_{B_{i_2}} 
\overline{\mathcal{M}}(B_{i_2},B_{i_3}) \times_{B_{i_3}} \cdots \times_{B_{i_{n-1}}}
\overline{\mathcal{M}}(B_{i_{n-1}},B_{i_n})
$$
where $m \geq i_1 > i_2 > \cdots > i_n \geq 0$, $Q$ is a face of $I^N$ of
dimension $q\leq p$, $\sigma:Q \rightarrow B_{i_1}$ is smooth, and the
fibered products are taken with respect to $\sigma$ and the beginning and
endpoint maps $\partial_-$ and $\partial_+$.   Lemma 5.1 of \cite{BanMor} shows
that the elements of $C_p$ are all compact smooth manifolds with corners. Let $S_p$
be the free abelian group generated by the elements of $C_p$, and let $S_p^\infty(B_i)$
denote the subgroup of the singular  $C_p$-chain group $S_p(B_i)$ generated by those
maps $\sigma:P \rightarrow B_i$ that satisfy the following two conditions: 
\begin{enumerate}
\item The map $\sigma$ is smooth.
\item If $P\in C_p$ is a connected component of a fibered product, 
      then $\sigma = \partial_+ \circ \pi$, where $\pi$ denotes projection 
      onto the last component of the fibered product.
\end{enumerate}

\begin{definition}\label{MBdegree}
Define the \textbf{Morse-Bott degree} of the singular topological chains
in $S_p^\infty(B_i)$ to be $p+i$. For any $k=0,\ldots ,m$ the group 
of smooth singular topological chains of Morse-Bott degree $k$ is defined to be
$$
\tilde{C}_k(f) = \bigoplus_{i=0}^m S_{k-i}^\infty(B_i).
$$
\end{definition}

If $\sigma:P \rightarrow B_i$ is a singular $C_p$-space in $S_p^\infty(B_i)$,
then for any $j = 1,\ldots ,i$ composing the projection map $\pi_2$ onto the second 
component of $P \times_{B_i} \overline{\mathcal{M}}(B_i,B_{i-j})$ with the 
endpoint map $\partial_+:\overline{\mathcal{M}}(B_i,B_{i-j}) \rightarrow B_{i-j}$ 
gives a map
$$
P \times_{B_i} \overline{\mathcal{M}}(B_i,B_{i-j}) 
\stackrel{\pi_2}{\longrightarrow} \overline{\mathcal{M}}(B_i,B_{i-j})
\stackrel{\partial_+}{\longrightarrow} B_{i-j}.
$$
Lemma 5.3 of \cite{BanMor} shows that restricting this map to the connected components
of the fibered product $P \times_{B_i}\overline{\mathcal{M}} (B_i,B_{i-j})$ and 
adding these restrictions (with the sign determined by the orientation 
when the dimension of a component is zero) defines an element 
$\partial_j(\sigma) \in S_{p+j-1}^\infty(B_{i-j})$.

\begin{definition}\label{boundary}
For $k=1,\ldots ,m$ define a homomorphism $\partial:\tilde{C}_k(f) 
\rightarrow \tilde{C}_{k-1}(f)$ as follows. If $\sigma \in S_{p}^\infty(B_i)$
is a singular $S_{p}$-space of $B_i$ where $p=k-i$, then
$$
\partial(\sigma) = \bigoplus_{j=0}^m \partial_j(\sigma)
$$
where $\partial_0$ is $(-1)^{k}$ times the boundary operator on singular
topological chains defined above, $\partial_j(\sigma) = \partial_+ \circ \pi_2:
P \times_{B_i} \overline{\mathcal{M}}(B_i,B_{i-j}) \rightarrow B_{i-j}$ for 
$j=1,\ldots ,i$, and $\partial_j(\sigma) = 0$ otherwise. The map $\partial$ extends
to a homomorphism
$$
\partial:\bigoplus_{i=0}^m S_{k-i}^\infty(B_i) \longrightarrow
\bigoplus_{i=0}^m S_{k-1-i}^\infty(B_i).
$$
\end{definition}

\noindent
The following is Proposition 5.5 of \cite{BanMor}.

\begin{proposition}\label{complex}
For every $j = 0, \ldots ,m$ we have $\sum_{q=0}^j \partial_q \partial_{j-q} = 0$.
\end{proposition}

\smallskip
Defining the Morse-Bott-Smale multicomplex over the integers requires a coherent
system of orientations on the elements of $C_p$ (cf. Section 5.2 of \cite{BanMor}) and
a collection of degeneracy relations that identify maps from different domains that are
``essentially'' the same (cf. Section 5.3 of \cite{BanMor}). The degeneracy relations 
are expressed in the form of subgroups $D^\infty_p(B_i)\subseteq S^\infty_p(B_i)$, and
the chain groups that make up the Morse-Bott-Smale multicomplex are defined to be
$S^\infty_p(B_i)/D^\infty_p(B_i)$. Lemma 5.10 of \cite{BanMor} shows that the homomorphisms
$\partial_j$ on $S^\infty_p(B_i)$ induce homomorphisms on $S_p^\infty(B_i) /D_p^\infty(B_i)$,
which we denote using the same notation.

\begin{definition}\label{MBScomplex}
Define
$$
C_p(B_i) = S_p^\infty(B_i) / D_p^\infty(B_i)
$$
to be the group of \textbf{non-degenerate} smooth singular topological chains in
$S_p^\infty(B_i)$. The group $C_k(f)$ of \textbf{$k$-chains} in the Morse-Bott
chain complex of $f$ is defined to be the group of non-degenerate smooth singular
topological chains of Morse-Bott degree $k$, i.e.
$$
C_k(f) = \bigoplus_{i=0}^m C_{k-i}(B_i) = \bigoplus_{i=0}^m
S_{k-i}^\infty(B_i)/D_{k-i}^\infty(B_i).
$$
The boundary operator in the Morse-Bott-Smale chain complex
$$
\partial:\bigoplus_{i=0}^m S_{k-i}^\infty (B_i)/D_{k-i}^\infty(B_i) \longrightarrow
\bigoplus_{i=0}^m S_{k-1-i}^\infty(B_i)/D_{k-1-i}^\infty(B_i)
$$
is defined to be $\partial = \oplus_{j=0}^m \partial_j$.
\end{definition}

\begin{figure}[h]
The Morse-Bott-Smale Multicomplex
$$
\xymatrix{
\ddots & \vdots\\
\cdots & C_1(B_2) \ar@{}[u]|{\oplus} \ar[r]^{\partial_0} \ar[dr]^{\partial_1} \ar[ddr]|(.33){\partial_2} |!{[d];[dr]}\hole & C_0(B_2) \ar[r]^{\partial_0} \ar[dr]^{\partial_1} \ar[ddr]|(.33){\partial_2} |!{[d];[dr]}\hole & 0 & & \\
\cdots & C_2(B_1)\ar@{}[u]|{\oplus} \ar[r]^(.55){\partial_0} \ar[dr]^{\partial_1} & C_1(B_1)\ar@{}[u]|{\oplus} \ar[r]^(.55){\partial_0} \ar[dr]^{\partial_1} & C_0(B_1) \ar@{}[u]|{\oplus} \ar[r]^{\partial_0} \ar[dr]^{\partial_1} & 0 & \\
\cdots & C_3(B_0) \ar@{}[u]|{\oplus} \ar[r]^{\partial_0} & C_2(B_0) \ar@{}[u]|{\oplus} \ar[r]^{\partial_0} & C_1(B_0) \ar@{}[u]|{\oplus} \ar[r]^{\partial_0} & C_0(B_0) \ar@{}[u]|{\oplus}  
 \ar[r]^-{\partial_0} & 0\\
\cdots &  {C}_3(f) \ar@{}[u]|{\|} \ar[r]^{\partial} & {C}_2(f) \ar@{}[u]|{\|} \ar[r]^{\partial} & {C}_1(f) \ar@{}[u]|{\|} \ar[r]^{\partial} & {C}_0(f) \ar@{}[u]|{\|} \ar[r]^-{\partial} & 0
}
$$
\end{figure}

\noindent
Since the homomorphisms $\partial_j$ are induced from the homomorphisms in Definition
\ref{boundary}, Proposition \ref{complex} shows that the relations that define a 
multicomplex are satisfied.


\medskip
\subsection*{Maps between Morse-Bott-Smale multicomplexes}
The Banyaga-Hurtubise approach to constructing the Morse-Bott-Smale multicomplex
has several advantages. For instance, it does not require picking any triangulations.
Hence, all the maps in the multicomplex are well defined at the chain level and
there is no need to prove that the homology of the multicomplex is independent of
arbitrarily chosen triangulations. Moreover, Lemma 5.1 of \cite{BanMor} shows that all
the fibered products used to define the multicomplex are compact smooth manifolds with
corners, without having to perturb any maps used in the construction. Other approaches
require perturbing the beginning and endpoint maps $\partial_-$ and $\partial_+$, which
would then necessitate proving that the homology of the resulting complex is independent
of the chosen perturbations \cite{FukFlo}.

\smallskip
While the multicomplex constructed by Banyaga and Hurtubise does not depend on any
extraneous choices, it obviously does depend on the Morse-Bott-Smale function
$f:M \rightarrow \mathbb{R}$ and the Riemannian metric on $M$. However, Theorem 6.17 of
\cite{BanMor} shows that the homology of the assembled chain complex does not depend on
the Morse-Bott-Smale function or the Riemannian metric on $M$. The proof of Theorem 6.17
of \cite{BanMor} follows standard continuation arguments found in papers on Floer homology. 
In particular, given two Morse-Bott-Smale functions $f_1$ and $f_2$ on $M$ a continuation
map is defined between the multicomplexes determined by the two functions using moduli 
spaces of time dependent gradient flow lines, i.e. moduli spaces of gradient flow lines
of a function $F_{21}:M \times \mathbb{R} \rightarrow
\mathbb{R}$ where 
\begin{eqnarray*}
\lim_{t \rightarrow -\infty} F_{21}(x,t) & = & f_1(x) + 1\\
\lim_{t \rightarrow +\infty} F_{21}(x,t) & = & f_2(x) - 1
\end{eqnarray*}
for all $x \in M$.

However, the time dependent moduli spaces of gradient flow lines are not allowed domains
for the singular topological chains in Morse-Bott-Smale multicomplex. So, Banyaga and
Hurtubise adapt the technique of \textbf{representing chain systems} from \cite{BarLag}
in order to define their continuation maps. Roughly speaking, a representing chain system
consists of singular topological chains (defined on the allowed domains) that represent the
fundamental classes of the moduli spaces of time dependent gradient flow lines (which are
compact smooth manifolds with corners). This means that the continuation maps are only defined
at the chain level after choosing a representing chain system. However, Corollary 6.12 
of \cite{BanMor} shows that the induced map between the homologies of the assembled chain
complexes is independent of the representing chain systems. So, the continuation maps are
well defined at the level of homology and independent of any of the choices made to define
them at the chain level.

\smallskip
The following two corollaries proved in Section 6 of \cite{BanMor} show that standard
arguments from Floer homology can be applied to the Morse-Bott-Smale multicomplex.

\begin{corollary}\label{canonical}
For any two Morse-Bott-Smale functions $f_1,f_2:M \rightarrow \mathbb{R}$
the time-dependent gradient flow lines from $f_1$ to $f_2$ determine
a canonical homomorphism 
$$
(F_{21})_\ast:H_\ast(C_\ast(f_1),\partial) \rightarrow H_\ast(C_\ast(f_2),\partial), 
$$
i.e. the map $(F_{21})_\ast$ is independent of the choice of the function
$F_{21}:M \times \mathbb{R} \rightarrow \mathbb{R}$ and the representing
chain system used to define the chain map $(F_{21})_\Box:C_\ast(f_1)
\rightarrow C_\ast(f_2)$.
\end{corollary}

\begin{corollary}\label{functorial}
For any four Morse-Bott-Smale functions $f_k:M \rightarrow \mathbb{R}$, where
$k=1,2,3,4$, the canonical homomorphisms satisfy
$$
(F_{43})_\ast \circ (F_{31})_\ast = (F_{42})_\ast \circ (F_{21})_\ast
$$
and
$$
(F_{32})_\ast \circ (F_{21})_\ast =  (F_{31})_\ast.
$$
\end{corollary}

\noindent
The preceding two corollaries and the Singular $N$-Cube Chain Theorem
(Theorem \ref{cubehomology}) imply the following, which is Theorem 6.17 
of \cite{BanMor}.

\begin{theorem}\label{homologyindependence}
The homology of the Morse-Bott chain complex $(C_\ast(f),\partial)$ is
independent of the Morse-Bott-Smale function $f:M \rightarrow \mathbb{R}$.
Therefore, 
$$
H_\ast(C_\ast(f),\partial) \approx H_\ast(M;\mathbb{Z}).
$$
\end{theorem}

\subsection*{Interpolating between singular $N$-cube chains and Morse chains}
When the function $f:M \rightarrow \mathbb{R}$ is Morse-Smale the critical
set $B_i$ is a discrete set of points for all $i=0,\ldots,m$, and the
groups $C_p(B_i)$ are trivial for all $p > 0$. When the function is constant
the entire manifold $M$ is a critical submanifold of Morse-Bott index zero.
In this case $B_i = \emptyset$ for all $i > 0$, and the groups $C_p(B_i)$ 
are trivial for all $i > 0$. These two cases appear in the diagram of a general
Morse-Bott-Smale chain complex as follows.
$$
\xymatrix{
&  \ddots \ar[dr]^{\partial_1}& & & & \\
&  &  C_0(B_2) \ar[dr]^{\partial_1} &  & & \\
&  & & C_0(B_1) \ar[dr]^{\partial_1} &  & \\
\cdots \ar[r]^{\partial_0} & C_3(B_0) \ar[r]^{\partial_0} & C_2(B_0) 
\ar[r]^{\partial_0} & C_1(B_0)  \ar[r]^{\partial_0} & C_0(B_0) \ar[dr]^{\partial_1}
\ar[r]^-{\partial_0} & 0\\
& & & & & 0
}
$$
In the first case the homomorphism $\partial_1$ is the Morse-Smale-Witten
boundary operator, and in the second case we have the chain complex of
singular $N$-cube chains, which computes the singular homology of $M$ by Theorem
\ref{cubehomology}. Thus, the Morse-Bott-Smale multicomplex provides a means of
interpolating between the Morse-Smale-Witten chain complex and the chain complex of
singular $N$-cube chains. Moreover, Theorem \ref{homologyindependence} shows that the
homology of these two chain complexes are the same, and hence the results in
\cite{BanMor} give a new proof of the Morse Homology Theorem (Theorem \ref{Morsehomology}).

\bibliographystyle{amsplain}
\bibliography{books,papers}

\end{document}